\documentclass[a4paper,12pt]{article}

\usepackage{lmodern}
\usepackage[x11names]{xcolor}
\usepackage[english]{babel}
\usepackage[shortlabels]{enumitem}

\usepackage{latexsym}

\usepackage[x11names]{xcolor}
\usepackage[utf8]{inputenc}
\usepackage[T1]{fontenc}
\usepackage{csquotes}

\usepackage{amsmath}
\usepackage{amsfonts}
\usepackage{amsbsy}
\usepackage{amssymb}
\usepackage{amsthm}
\usepackage{amscd}
\usepackage{authblk}

\usepackage{bookmark} 

\usepackage{marvosym} \usepackage{mathtools} \usepackage{stmaryrd}

\usepackage{tabularx}
\usepackage[font=footnotesize]{caption}
\usepackage{subcaption}
\usepackage{colortbl}

\usepackage{tikz-cd}

\theoremstyle{plain}
\newtheorem{theorem}{Theorem}[section]

\theoremstyle{definition}

\newtheorem{assumption}[theorem]{Assumption}

\theoremstyle{remark}
\newtheorem{remark}[theorem]{Remark}

\numberwithin{equation}{section}

\allowdisplaybreaks

\newcommand{\dd}{\mathrm{d}}
\newcommand{\ii}{\mathrm{i}}

\newcommand{\Ac}{\mathcal A}

\newcommand{\Cc}{\mathcal C}

\newcommand{\Fc}{\mathcal F}

\newcommand{\Lc}{\mathcal L}

\newcommand{\Nc}{\mathcal N}
\newcommand{\Oc}{\mathcal O}

\newcommand{\Xc}{\mathcal X}

\newcommand{\C}{\mathbb C}

\newcommand{\N}{\mathbb N}

\newcommand{\R}{\mathbb R}

\newcommand{\Z}{\mathbb Z}

\newcommand{\Eb}{\mathbf{E}}
\newcommand{\Fb}{\mathbf{F}}

\newcommand{\Pb}{\mathbf{P}}
\newcommand{\Qb}{\mathbf{Q}}

\newcommand{\Sb}{\mathbf{S}}

\newcommand{\Yb}{\mathbf{Y}}

\usepackage{mathrsfs}

\newcommand{\As}{\mathscr A}
\newcommand{\Bs}{\mathscr B}

\newcommand{\Fs}{\mathscr F}

\newcommand{\Xs}{\mathscr X}

\newcommand{\Ef}{\mathfrak E}

\usepackage{dsfont}
\newcommand{\1}{\mathbf 1}

\renewcommand{\epsilon}{\varepsilon}
\renewcommand{\phi}{\varphi}

\newcommand{\KL}{\mathrm{KL}}

\newcommand{\Cov}{\operatorname{Cov}}

\newcommand{\defeq}{\vcentcolon=}

   \newcommand{\abs}[1]{\vert#1\vert}   

\usepackage{mathtools}

 \newcommand{\Expm}{\Ef}
\newcommand{\FctClass}{\Theta}

\newcommand{\gn}{\xi}

\newcommand{\Ytilde}{\widetilde Y}

\newcommand{\xb}{\boldsymbol x}
\newcommand{\yb}{\boldsymbol y}

\newcommand{\Xib}{\boldsymbol \Xi}

\newcommand{\RKHS}{\mathcal H}

\newcommand{\eb}{\boldsymbol{\mathrm e}}
\newcommand{\cdotb}{\boldsymbol{\cdot}}
\newcommand{\alphab}{{\boldsymbol \alpha}}

\usepackage[paper=a4paper,left=30mm,right=30mm,top=30mm,bottom=30mm]{geometry}

\usepackage{hyperref}
\hypersetup{colorlinks=true,urlcolor=RoyalBlue4,linkcolor=RoyalBlue4,citecolor=RoyalBlue4,naturalnames=true,hypertexnames=false}

\usepackage{titlesec}

\titleformat*{\section}{\large\bfseries}
\titleformat*{\subsection}{\large}
\titleformat*{\subsubsection}{\large\bfseries}
\titleformat*{\paragraph}{\large\bfseries}
\titleformat*{\subparagraph}{\large\bfseries}

\usepackage[backend=bibtex,giveninits, style=alphabetic, uniquelist=false, natbib=true, doi=false, url=false, isbn = false]{biblatex}

\DeclareFieldFormat[article]{title}{#1} \DeclareFieldFormat[book]{title}{\itshape #1}
\DeclareFieldFormat[thesis]{title}{#1}
\DeclareFieldFormat[incollection]{title}{\itshape #1}
\DeclareFieldFormat[inproceedings]{title}{\itshape #1}
\DeclareFieldFormat[article]{author}{\color{blue} #1}
\DeclareFieldFormat[article]{pages}{#1.}
\DeclareFieldFormat[incollection]{pages}{#1.} 

\renewbibmacro{in:}{\ifentrytype{article}{}{\printtext{\bibstring{in}\intitlepunct}}}

\DeclareNameAlias{sortname}{last-first}
\DeclareNameAlias{default}{last-first}

\renewbibmacro*{volume+number+eid}{\printfield{volume}\setunit*{\addnbthinspace}\printfield{number}\setunit{\addcomma\space}\printfield{eid}}
\DeclareFieldFormat[article]{number}{\mkbibparens{#1}}

\renewbibmacro*{publisher+location+date}{\printlist{publisher}\iflistundef{location}
	{\setunit*{\addcomma\space}}
	{\setunit*{\addcomma\space}}\printlist{location}\setunit*{\addcomma\space}\usebibmacro{date}\newunit}

 \DeclareFieldFormat[article]{number}{}

\title{Asymptotic equivalence for nonparametric regression with dependent errors: Gauss-Markov processes\thanks{This research has been supported in part by the research grant  DE 502/27-1 of the German Research Foundation (DFG).}}
\author{Holger Dette}
\author{Martin Kroll}
\affil{Ruhr-Universität Bochum}
 
\date{\today}

\begin{document}

\maketitle

\begin{abstract}
  For the   class of  Gauss-Markov processes we  study the  problem of asymptotic equivalence of 
  the nonparametric regression model with errors given by the increments of the  process and the continuous time model, where
  a whole path of a  sum of a deterministic signal and the  Gauss-Markov process can be observed.
 In particular we provide  sufficient conditions such that asymptotic equivalence of the two models holds for functions from a given class, and we verify these  for the special cases of Sobolev ellipsoids and Hölder classes with smoothness index $> 1/2$ under mild assumptions on the Gauss-Markov process at hand.
    To derive these results, we develop an explicit characterization of the reproducing kernel Hilbert space associated with the  Gauss-Markov process, that hinges on a characterization of such processes by a property of the corresponding covariance kernel
    introduced by Doob \citep{doob1949heuristic}.
    In order to demonstrate that the given assumptions on the Gauss-Markov process are in some sense sharp we also
    show that asymptotic equivalence fails to hold for the special case of Brownian bridge.
    Our findings demonstrate that the well-known asymptotic equivalence of the Gaussian white noise model and the nonparametric regression model with i.i.d.\,standard normal errors (see Brown and Low, \citep{brown1996asymptotic}) can be extended to a result treating general Gauss-Markov noises in a unified manner.
\end{abstract}



\normalsize

\section{Introduction}

In a seminal paper Brown and Low \citep{brown1996asymptotic} establish asymptotic equivalence of  the nonparametric  regression model with discrete observations
\begin{equation}\label{eq:BL:regr}
	Y_{i,n} = f \left( \frac{i}{n} \right) + \eta_i, \qquad i = 1,\ldots,n,
\end{equation}
and the continuous time model defined by the stochastic differential equation
\begin{equation}\label{eq:BL:GWN}
	\dd Y_t = f(t)\dd t + \frac{1}{\sqrt{n}} \dd W_t, \quad t \in [0,1], 
\end{equation}
 where $\eta_1, \ldots , \eta_n$ are independent, standard Gaussian random variables, 
  $W$ denotes a standard Brownian motion (that is, $\dd W$ is white noise) and $f$ is the unknown nonparametric drift 
 satisfying a smoothness assumption.
 Equation~\eqref{eq:BL:GWN} is commonly referred to as the \emph{Gaussian white noise model} and serves as an important benchmark model in nonparametric statistics.
Often, due to the absence of discretization effects, statistical methods are easier to analyze in model \eqref{eq:BL:GWN} than in \eqref{eq:BL:regr}.
Asymptotic equivalence between the two models then suggests that theoretical results obtained in the Gaussian white noise model hold true in the more realistic model \eqref{eq:BL:regr} as well.

Since the contribution of Brown and Low
numerous authors have worked  on the  problem of  establishing  asymptotic equivalence of various models from different perspectives. For example, Grama and Nussbaum \citep{gramanussbaum1998} investigate nonparametric generalized linear models, and Brown and Zhang \citep{brown1998asymptotic} prove nonequivalence when the smoothness of the function class is equal to $ 1/2$.  Brown et al.\,\citep{brown2002asymptotic} and Rei\ss \, \citep{reiss2008asymptotic}
study  the random  design in the one-dimensional and multivariate case, respectively 
(see also Carter, \citep{carter2006} for some results in models with a multivariate fixed design).
The general framework in \citep{brown1996asymptotic} is already formulated for heteroscedastic errors and Carter \citep{carter2007asymptotic} shows
asymptotic equivalence for unknown variances and design density. We also refer to the work of  
Rei\ss \, \citep{reiss2011asymptotic} and Meister \citep{meister2011asymptotic}
who propose rate-optimal estimators of the volatility function and sharp minimax constants in the functional linear regression model as an application of asymptotic equivalence, respectively.

This list is by no means complete, but a common feature of most publications in this field consists in the fact that the random variables in the corresponding discrete nonparametric regression model are assumed to be independent. Grama and Neumann \citep{GramaNeumann2006} consider a nonparametric autoregression model but still use an i.i.d.\,assumption for the innovations.
Golubev, Nussbaum, and Zhou \citep{golubev2010asymptotic} study a stationary process and show asymptotic equivalence in the context of spectral density estimation.
 However, when it comes to regression the assumption of i.i.d.\,errors is often made for the theoretical analysis.
 Motivated by the work of Johnstone and Silverman \citep{johnstone1997wavelet} and Johnstone \citep{johnstone1999wavelet}, Carter \citep{carter2009asymptotically} considers asymptotic equivalence under the  assumption that the noise process is given by a wavelet composition with stochastically independent coefficients.
 Another notable exception is the recent work of Schmidt-Hieber \citep{schmidt-hieber2014asymptotic}, who considers the  
nonparametric regression model  \eqref{eq:BL:regr} with fractional Gaussian noise (fGN) and establishes the asymptotic equivalence to a model of the form \eqref{eq:BL:GWN}, where the error process $\frac{1}{\sqrt n} W_t$ is replaced by a fractional Brownian motion (fBM)  $n^{H-1} B_t^H$ with Hurst parameter 
$H \in (\frac{1}{4},1)$ over periodic Sobolev ellipsoids containing sufficiently smooth functions.
It is also shown  that asymptotic equivalence fails to hold for certain combinations of Hurst parameter and smoothness index leading to sharp results concerning the smoothness requirement in the case $H \in [\frac{1}{2},1)$.

The purpose of the present paper is to provide an essentially distinct way of investigating  
nonparametric regression models with dependent errors for asymptotic equivalence. 
Instead of replacing the  Brownian motion by a  fractional Brownian motion, we consider here arbitrary Gauss-Markov processes as error processes in model \eqref{eq:BL:GWN}.
To be  precise let $\Xi = (\Xi_t)_{t \in [0,1]}$ be such a Gauss-Markov process with initial state $\Xi_0 = 0$, and consider the continuous time model
\begin{equation}\label{eq:intro:cont}
	Y_t = \int_0^t f(s) \dd s + \frac{1}{\sqrt{n}} \Xi_t, \quad t \in [0,1].
\end{equation} 
Assuming equidistant design points, the candidate regression model for asymptotic equivalence is given by 
\begin{equation}\label{eq:intro:regr}
		Y_{i,n} = f  \left( \frac{i}{n} \right)  +  \sqrt{n} \xi_{i,n}, \quad i=1,\ldots,n,
\end{equation}
where $\xi_{i,n} \stackrel{\Lc}{=} \Xi_{i/n}  - \Xi_{(i-1)/n}$ denote the increments of the process $ \Xi$. Note that in general  the observation errors $\xi_{i,n}$ are not uncorrelated. However, for the special case of Brownian motion, that is $\Xi   \stackrel{\Lc}{=} W$, Equations \eqref{eq:BL:regr} and \eqref{eq:intro:regr} are equivalent in distribution since Brownian motion has independent increments and satisfies the scaling property that  $W(\cdotb)$ and $c^{-1/2}W(c\,\cdotb)$ have the same distribution for any $c>0$.
The main results of the present  paper establish the asymptotic equivalence of the models  \eqref{eq:intro:cont} and \eqref{eq:intro:regr} for a wide class of Gauss-Markov processes for functional parameters $f$ belonging to a Sobolev or Hölder class of sufficiently smooth functions.

There are different ways to  prove asymptotic equivalence between regression and white noise experiments in the literature.
The original paper of Brown and Low \citep{brown1996asymptotic}
considers the case where the regression $f$ is an element of a function class, say  $\Theta$, and
uses the  key assumption 
\begin{equation}\label{eq:BL:cond}
	\lim\limits_{n \to \infty} \sup_{f \in \Theta} n \int_0^1 \Big (f(t) - \bar{f}_n(t) \Big )^2\dd t = 0,
\end{equation}
where $\bar f_n(\cdotb) = \sum_{j=1}^{n} f(j/n) \1_{((j-1)/n,j/n]}(\cdotb)$  denotes a piecewise constant approximation of $f$ and 
$\1_A$ is the indicator function of the set $A$.
These  authors  introduce an intermediate set of random variables that forms a sufficient statistic for the white noise model with the function $f$ replaced with $\bar f_n$.
Then, the Hellinger distance between this sufficient statistic and the regression experiment with discrete sampling locations is shown to converge to zero.
On the other hand, condition \eqref{eq:BL:cond} guarantees that the Le Cam distance between the white noise models with parameters $f$ and $\bar f_n$, respectively, vanishes as $n \to \infty$.
Note that this approach  can be used to  show asymptotic equivalence but cannot provide the optimal rate of convergence for the Le Cam distance between the two experiments.
This issue has been solved  by Rohde \citep{rohde2004asymptotic}, who considers a Gaussian sequence space model as an intermediate experiment between the two experiments of interest. 
Schmidt-Hieber \citep{schmidt-hieber2014asymptotic} in some sense generalizes the conditions in \citep{brown1996asymptotic} and formalizes them in the framework of reproducing kernel Hilbert spaces (RKHSs), which is a suitable setup for the investigation in the case of fractional Brownian motion.
Our analysis of the asymptotic equivalence for models with Gauss-Markov errors
will also be based on the RKHS framework.  As an essential ingredient  we will use  a characterization  of Gauss-Markov processes 
introduced by Doob \citep{doob1949heuristic} to 
derive  an explicit representation of the RKHSs associated with the Gauss-Markov processes under consideration, which 
can be used to develop  sufficient conditions  for the  asymptotic equivalence of the models \eqref{eq:intro:cont} and \eqref{eq:intro:regr}.

The remaining part of the paper is organized as follows.
In Section~\ref{s:GM} we introduce Gauss-Markov processes and recap the  characterization of such processes introduced by Doob \citep{doob1949heuristic} that will be pivotal for our approach. Roughly speaking, these processes can be characterized by the property that the corresponding covariance kernel is {\it  triangular}  (see Section \ref{subs:GM} for a precise definition).
In Section~\ref{s:RKHS} we study the reproducing kernel Hilbert spaces (RKHSs) associated with Gauss Markov processes and derive representations of these spaces via Hilbert space isomorphisms into a space of square-integrable functions.
In Section~\ref{s:cond} we recall the basic elements  of Le Cam theory and provide sufficient conditions on the class of potential functions $f$ and the Gauss-Markov process that imply asymptotic equivalence.
The characterizations of the RKHSs are used in Section~\ref{s:sobolev} where we establish asymptotic equivalence of  the models 
\eqref{eq:intro:cont} and \eqref{eq:intro:regr} under mild assumptions on the Gauss-Markov process in model \eqref{eq:intro:cont} for Sobolev ellipsoids and H\"older classes.
In Section~\ref{subs:sketch} we discuss a different approach to derive such results, which is not based 
on RKHS theory. It requires a (slightly) stronger assumption than those used in this paper and 
applies results from the seminal paper \citep{brown1996asymptotic} to suitable transformations of the experiments of interest.

Finally, in Section~\ref{s:BB} we demonstrate  that asymptotic equivalence cannot hold without any additional assumptions on the Gauss-Markov process.
More precisely, for the special case of Brownian bridge we show that the Le Cam distance between the two experiments is bounded away from zero.
The proofs of our results are deferred to the Appendix.
\subsection*{Notation} 
Vectors will be denoted with bold letters (i.e., we write $\xb_n = (x_{1,n},\ldots,x_{m(n),n})$ when both length and entries of a vector might vary with $n$).
We also use the shorthand $\xb_{1:k,n} = (x_{1,n},\ldots,x_{k,n})$.
We write $a_n \lesssim b_n$ if $a_n \leq C b_n$ for a constant $C$ that is independent of $n$.
The shorthand $a_n \asymp b_n$ is used when $a_n \lesssim b_n$ and $b_n \lesssim a_n$ hold simultaneously. 
\section{The RKHS associated with a Gauss-Markov process}\label{s:GM}

The purpose of this section is to lay the foundations for our main results in Sections~\ref{s:cond} and \ref{s:sobolev}.
In Section~\ref{subs:GM} we recapitulate a characterization of Gauss-Markov process going back to Doob which will be important for our further reasoning.
As a consequence of this characterization we can give an explicit description of the RKHS associated with the covariance kernel of a Gauss-Markov process, which is of independent interest and presented in Section \ref{s:RKHS}.

\subsection{Gauss-Markov processes}\label{subs:GM}

By definition a Gauss-Markov process is a stochastic process that is both Gaussian and Markov.
Such a process $X = (X_t)_{t \in [0,1]}$ is essentially characterized by the following factorization property of the covariance function:
\begin{equation}\label{eq:tri:kernel}
K_X(s,t) \defeq \Cov(X_s,X_t) = \Eb [X_s X_t] = U(s)V(t) \quad \text{for } 0 \leq s \leq t \leq 1,
\end{equation}
where $U$ and $V$ are (known) non-negative functions on the interval $[0,1]$; see \cite{doob1949heuristic,mehr1965certain} for details.
Kernels with the factorization property \eqref{eq:tri:kernel} are sometimes referred to as \emph{triangular kernels} in the literature.
Examples of Gauss-Markov processes include standard Brownian motion ($U(t) = t$, $V \equiv 1$), the Ornstein-Uhlenbeck process ($U(t)=\exp(Lt)$ and $V(t) = \exp(-Lt)$ for some $L > 0$) and the Brownian bridge ($U(t)=t$, $V(t) = 1-t$).

For our results, we will further assume that the considered Gauss-Markov process starts in zero.
Such a process will be denoted with $\Xi= (\Xi_t)_{t \in [0,1]}$ from now on, and we impose the following assumption.
\begin{assumption}\label{ass:GM}
	The process $\Xi = (\Xi_t)_{t \in [0,1]}$ is a Gauss-Markov process with $\Xi_0 = 0$, and non-degenerate on the open interval $(0,1)$.
\end{assumption}
Note that Assumption~\ref{ass:GM} implies that there exist functions $u$ and $v$ in the representation
\begin{equation}\label{eq:tri:kernel:0}
    K_\Xi(s,t)=u(s)v(t)
\end{equation}
of the covariance kernel satisfying $u(\cdotb)v(\cdotb) \geq 0$ on the interval $[0,1]$, $u(\cdotb)v(\cdotb) > 0$ on $(0,1)$ and that  the function
	\begin{equation*}
		q(t) = \frac{u(t)}{v(t)},
	\end{equation*}
	is continuous on   the interval  $[0,1)$, non-negative and strictly increasing on $[0,1]$ (see \cite{mehr1965certain}, p.~507).
Moreover, under Assumption~\ref{ass:GM}, the Gauss-Markov process $\Xi$ can be written in distribution as
\begin{equation}\label{eq:repr:GM}
\Xi_t = v(t) \cdot W_{q(t)}.
\end{equation}
Vice versa, this transformation of Brownian motion defines a Gaussian process with covariance kernel \eqref{eq:tri:kernel:0}.
Note that a simple calculation shows that the process $\Xi$ has independent increments if and only if the function $v$ is constant.

Let us shortly explain how processes $\Xi$ satisfying Assumption~\ref{ass:GM} can be obtained.
First, one can of course define such a process directly by means of the representation \eqref{eq:repr:GM} provided that the function $u$ and $v$ satisfy the properties stated above and 
$u(0) = 0$ or $v(0)=0$ (this latter assumption guarantees that the process starts in $0$).
Second, starting with a general centered non-degenerate Gauss-Markov process $X=(X_t)_{t \in [0,1]}$ with covariance kernel \eqref{eq:tri:kernel}, and then conditioning on $X_0 = 0$ does neither suspend Gaussianity nor the Markov property.
More precisely, the process $\Xi = (\Xi_t)_{t \geq 0}$ with $\Xi_t \sim X_t| \{ X_0=0 \}$ is a centered Gaussian process with covariance kernel
\begin{equation*}
	K_\Xi (s,t) \defeq \Cov(\Xi_s,\Xi_t) = \Eb [\Xi_s \Xi_t] = u(s)v(t), \qquad 0 \leq s \leq t \leq 1,
\end{equation*}
where the functions $u$ and $v$ are related to $U$ and $V$ through the identities
\begin{equation*}
	u(t) = U(t) - Q(0)V(t) \qquad \text{and} \qquad v(t) = V(t),
\end{equation*}
where $Q(t)=U(t)/V(t)$.
In particular, the process $\Xi$ satisfies Assumption~\ref{ass:GM}, and $q(t)=u(t)/v(t)$ can be represented as $q(t) = Q(t) - Q(0)$.

\medskip

Surprisingly, Gauss-Markov processes starting at zero can  also be obtained by conditioning a non-Markovian Gaussian process   $Y= (Y_t)_{t\ge 0} $ on the event $\{ Y_0=0 \}$.
Then, the initial Gaussian process $Y$ is called \emph{conditionally Markov}.
As one interesting example (further examples of conditionally Markov processes can again be found in \cite{mehr1965certain}, pp.~513 and 516) let us state the stationary Gaussian process $(Y_t)_{t \in \R}$ defined on the whole real-line with zero mean and covariance kernel
\begin{align*}
	K_Y(x,y) = \begin{cases}
		1 - \lvert x - y\rvert, & \text{ if } 0 \leq \lvert x-y\rvert \leq 1,\\
		0, & \text{ if } \lvert x-y\rvert > 1.
	\end{cases}
\end{align*}
for $x,y \in \R$.
Then, the process $(\Xi_t)_{t \geq 0}$ with $\Xi_t \sim Y_t|\{ Y_0 = 0 \}$ is a centered Gaussian process with covariance kernel
\begin{equation*}
	K_\Xi(t_1,t_2) = \begin{cases}
		t_1(2-t_2), & \text{ if } 0 \leq t_1 \leq t_2 \leq 1,\\
		1 - (t_2-t_1), & \text{ if }  0 < t_2 - t_1 < 1, t_2 \geq 1, t_1 \geq 0,\\
		0, & \text{ if } t_2 - t_1 \geq 0, t_1 \geq 0.
	\end{cases} 
\end{equation*}
In particular, the restriction of $\Xi$ on the interval $[0,1]$ is a Gauss-Markov process satisfying Assumption~\ref{ass:GM}, and $K_\Xi(s,t)=u(s)v(t)$ for $0 \leq s \leq t \leq 1$ with $u(x)=x$ and $v(x)=2-x$.
This process has been considered by Slepian in \cite{slepian1961first} where it was proved that the process fulfills a 'peculiar Markov-like property'.
The restriction of the process $\Xi$ on the unit interval fulfills all the technical assumptions made on the functions $u$ and $v$ below, and thus the results on  asymptotic equivalence  derived in Section \ref{s:cond} and \ref{s:sobolev} hold in particular for this process.

\subsection{Reproducing kernel Hilbert spaces}\label{s:RKHS}

We briefly recall some basic facts about reproducing kernel Hilbert spaces. For a detailed discussion the reader 
is referred to the monographs of Berlinet and Thomas-Agnan \cite{berlinet2004reproducing} and Paulsen and Raghupathi \cite{paulsen2016introduction}.
A subset $\RKHS$ of the set $\Fc(\Xc,\R)$ of real-valued functions on a domain $\Xc$ is called a reproducing kernel Hilbert space (RKHS) on $\Xc$ if
\begin{itemize}
	\item $\RKHS$ is a vector subspace of $\Fc(\Xc,\R)$,
	\item $\RKHS$ has an inner product $\langle \cdotb,\cdotb \rangle_{\RKHS}$ with respect to which $\RKHS$ is a Hilbert space,
	\item for any $x \in \Xc$, the linear evaluation functional $\ell_x \colon \RKHS \to \R$ with $\ell_x(f) = f(x)$ is bounded.
\end{itemize}
The Riesz representation theorem implies that for any $x \in \Xc$ there exists a unique function $k_x \in \RKHS$ such that for every $f \in \RKHS$ one has $f(x) = \ell_x(f) = \langle f,k_x \rangle_{\RKHS}$.
The function $K \colon \Xc \times \Xc \to \R$
\begin{equation*}
	K(x,y) = k_x(y)
\end{equation*}
is called the \emph{reproducing kernel} for $\RKHS$.
Reproducing kernels are \emph{positive definite} in the sense that the inequality
\begin{equation*}
	\sum_{i,j=1}^{n} a_i \overline a_j K(x_i,x_j) \geq 0
\end{equation*}
holds for all $n \in \N$, $x_i \in \Xc$ and $a_i \in \C$.
Vice versa, for any positive definite function $K$ on $\Xc$ there exists a RKHS $\RKHS(K)$ on $\Xc$ with reproducing kernel $K$ and this RKHS is uniquely determined by the kernel.
In particular, this holds true for covariance kernels.
Finally, it is well-known that the linear span of the functions $k_y=K(\cdotb,y)$ is dense in $\RKHS$.

In the following, we will consider the RKHS $\RKHS(\Xi)$ associated with the covariance kernel $K_\Xi\colon [0,1] \times [0,1] \to \R$ of a Gauss-Markov process $\Xi$ satisfying Assumption~\ref{ass:GM} as introduced in Section~\ref{s:GM}.
As a consequence
\begin{equation*}K_\Xi(s,t) = u(s)v(t), \qquad 0 \leq s \leq t \leq 1.
\end{equation*}
Let us from now on assume that $q(0)=0$ for $q(x)=u(x)/v(x)$ which is under the validity of Assumption~\ref{ass:GM} certainly satisfied whenever $v(0) \neq 0$.
This additional condition is valid for all the results in Section~\ref{s:sobolev}.
As mentioned above, $\Xi$ can be represented in distribution as
\begin{equation*}
	\Xi_t =  v(t)W_{q(t)}
\end{equation*}
where $W$ is  standard Brownian motion.
Denote $T=q(1)$, and let us consider the mapping $\psi \colon \RKHS(\Xi) \to L^2([0,T])$ defined on the 'generators' $ \{ K_\Xi(\cdotb, t)$, $t \in [0,1] \} $ via
\begin{equation} \label{psi} 
	K_\Xi(\cdotb, t) \stackrel{\psi}{\mapsto} v(t) \1_{[0,q(t)]}(\cdotb)
\end{equation}
(on general elements of $\RKHS(\Xi)$ is naturally defined by a limiting process). The mapping 
$\psi$ is an isometry of Hilbert spaces since
\begin{align*}
	\langle K_\Xi(\cdotb,s), K_\Xi(\cdotb,t) \rangle_{\RKHS(\Xi)} &= K_\Xi(s,t)
	= v(s)v(t)q(s)\\
	&= \int_0^T v(s)v(t) \1_{[0,q(s)]}(u) \1_{[0,q(t)]}(u) \dd u \\
	&= \langle \psi K_\Xi(\cdotb,s), \psi K_\Xi(\cdotb,t), \rangle_{L^2([0,T])}
\end{align*}
for $0 \leq s \leq t \leq 1$ (this holds true due to the strict monotonicity and continuity of $q$ combined with the fact that the indicator functions $\1_{[0,t]}, t \in [0,T]$, are dense in $L^2([0,T])$).
We use the shorthand notation $\psi_t = \psi(K_\Xi(\boldsymbol \cdotb,t))$.
For $F \in \RKHS(\Xi)$ we obtain with $g = \psi F \in L^2([0,T])$ the representation 
\begin{align*}
	F(t) &= \langle F, K_\Xi(\boldsymbol \cdot,t) \rangle_{\RKHS(\Xi)}
	= \langle \psi F, \psi_t \rangle_{L^2([0,T])}
	= \langle g, \psi_t \rangle_{L^2([0,T])}\\
	&= \int_0^T g(u) v(t) 1_{[0,q(t)]}(u) \dd u
	= v(t)\int_0^{q(t)} g(u) \dd u~.
\end{align*}
This gives an explicit characterization of the RKHS $\RKHS(\Xi)$, namely
\begin{equation}\label{eq:RKHS:explicit}
\RKHS(\Xi) = \left\lbrace F: [0,1] \to \R ~\Big | ~ F(t) = v(t) \int_0^{q(t)} g(u)\dd u \text{ for some } g \in L^2([0,T]) \right\rbrace,
\end{equation}
and this space is equipped with the norm $\lVert F \rVert_{\RKHS(\Xi)} = \lVert g \rVert_{L^2([0,T])}$ for $F$ as in \eqref{eq:RKHS:explicit}.
Of course, in the special case of Brownian motion ($u(t)=t, v \equiv 1$) these calculations yield the well-known fact that the RKHS 
corresponding to the kernel of standard BM contains exactly the primitives (starting at $0$) of square-integrable functions.
If $F(\cdotb)=F_f(\cdotb)= \int_0^{\cdotb} f(s) \dd s$, then $f$ can be derived by differentiation:
\begin{align*}
	f(t) &= \frac{\partial}{\partial t} \left[ v(t) \int_0^{q(t)} g(u)\dd u \right]\\
	&= v'(t) \int_0^{q(t)} g(u) \dd u + v(t) g(q(t)) q'(t).
\end{align*}
Moreover, the function 
$g$ can be obtained from $F_f$ (and thus from $f$) by
\begin{align}\label{eq:g:F}
	g(q(t)) &= \left( \frac{F_f(t)}{v(t)} \right)' \cdot \frac{1}{q'(t)}. \end{align}
This relation between $f$ and $g$ will be exploited when proving asymptotic equivalence results for concrete function classes in Section~\ref{s:sobolev}. 
\section{Abstract conditions for asymptotic equivalence}\label{s:cond}

In this section, we derive sufficient conditions for an abstract function class $\Theta$ that guarantee asymptotic equivalence of nonparametric regression with dependent errors and a continuous time model with an additive error process that is Gauss-Markov.
In Section~\ref{subs:lecam} we gather the necessary notions from asymptotic equivalence theory.
In Section~\ref{subs:suff}, we state the main result that provides the announced sufficient conditions for asymptotic equivalence of models \eqref{eq:intro:cont} and \eqref{eq:intro:regr}.

\subsection{Asymptotic equivalence of experiments}\label{subs:lecam}

Le Cam's equivalence theory for statistical experiments has by now become a common tool in nonparametric statistics.
As one of its major appeals one might consider the fact that complex statistical models can be shown equivalent to simple and well-studied benchmark models, at least asymptotically when the amount of information contained in the data increases.
This is not only of interest in its own but has also turned out useful when proving optimality properties of estimation techniques.
For instance, Rei\ss \, \cite{reiss2011asymptotic}  proposes rate-optimal estimators of the volatility function and simple efficient estimators of the integrated volatility as an application of Le Cam theory.
Similarly, Meister \cite{meister2011asymptotic} derives sharp minimax constants in the functional linear regression model as an application of asymptotic equivalence of this model and an inverse problem in a white noise setup.
The interested reader will find comprehensive introductions into Le Cam theory in the monographs \cite{lecam1986asymptotic,lecam2000asymptotics,torgersen1991comparison} as well as in the recent survey paper \cite{mariucci2016lecam} by Mariucci, the latter focusing on nonparametric models that are also the topic of the present paper. 

There exist various equivalent ways to introduce the concept of equivalence of statistical experiments.
One approach in the general theory introduced by Le Cam is by means of the abstract concept of \emph{transitions}.
In the case of dominated statistical models with Polish sample spaces (which is exclusively considered in this article), however, this general notion essentially boils down to the class of \emph{Markov kernels} as has been shown in Proposition~9.2 in \cite{nussbaum1996asymptotic}.  
To be precise, let $\Expm_{i} = (\Xc_i,\Xs_i,(\Pb_{i,\theta})_{\theta \in \Theta})$ for $i \in \{1,2\}$ be two such dominated experiments with Polish parameter spaces $(\Xc_i,\Xs_i)$, and sharing the same parameter space $\Theta$.
Then, the \emph{deficiency} of $\Expm_1$ with respect to $\Expm_2$ is defined as the quantity
\begin{equation*}
	\delta(\Expm_1,\Expm_2) = \inf_{K} \sup_{\theta \in \Theta} \lVert K\Pb_{1,\theta} - \Pb_{2,\theta} \rVert_{\mathrm{TV}}~,
\end{equation*}
where $\lVert \Pb - \Qb \rVert_{\mathrm{TV}}$ denotes the total variation distance between probability measures $\Pb$ and $\Qb$, and the infimum is taken over all Markov kernels $K\colon \Xc_1 \times \Xs_2 \to [0,1]$.
$\delta(\Expm_{1},\Expm_{2})=0$ has the interpretation that the experiment $\Expm_{1}$ is more informative than $\Expm_{2}$. 

Building on the definition of deficiency, the \emph{Le Cam distance} between $\Expm_1$ and $\Expm_2$ is defined by symmetrization as
\begin{equation*}
	\Delta(\Expm_1,\Expm_2) = \max \{ \delta(\Expm_1,\Expm_2), \delta(\Expm_2,\Expm_1) \},
\end{equation*}
which provides a pseudo-metric on the space of all statistical models with common parameter space.
Two experiments $\Expm_1$ and $\Expm_2$ are called \emph{equivalent} if $\Delta(\Expm_1,\Expm_2) = 0$.
More generally, two sequences $(\Expm_{1,n})_{n \in \N}$ and $(\Expm_{2,n})_{n \in \N}$ are said to be \emph{asymptotically equivalent} if $\Delta(\Expm_{1,n}, \Expm_{2,n}) \to 0$.
In the latter case we write $\Expm_{1,n} \approx \Expm_{2,n}$.

\subsection{Sufficient conditions for asymptotic equivalence}\label{subs:suff}

We now rigorously define the two statistical experiments that  will be examined for asymptotic equivalence in this paper.
Let $\Theta$ denote a given class of functions. The first experiment with discrete observations  is given by 
\begin{equation*}
	\Expm_{1,n}(\Theta) = (\R^n, \Bs(\R^n), (\Pb_{1,n}^f)_{f \in \Theta})~,
\end{equation*}
where $\Pb_{1,n}^f$ denotes the distribution of  the vector
$\Yb_n = (Y_{1,n},\ldots,Y_{n,n})$ with components $Y_{i,n}$ defined by
\begin{equation}\label{eq:Yin}
Y_{i,n} = f  ( t_{i,n} ) +  \sqrt{n} \xi_{i,n}, \quad i=1,\ldots,n,
\end{equation}
where $t_{i,n} = i/n$ are the sampling locations, and $\xi_{i,n} = \Xi \left( t_{i,n} \right) - \Xi \left( t_{i-1,n} \right)$ are the increments of a centered Gauss-Markov process $\Xi = (\Xi_t)_{t \in [0,1]}$ with $\Xi_0  = 0$ as introduced in Section~\ref{s:GM}.

\begin{figure}
     \centering
     \begin{subfigure}[b]{0.49\textwidth}
         \centering
         \includegraphics[width=\textwidth]{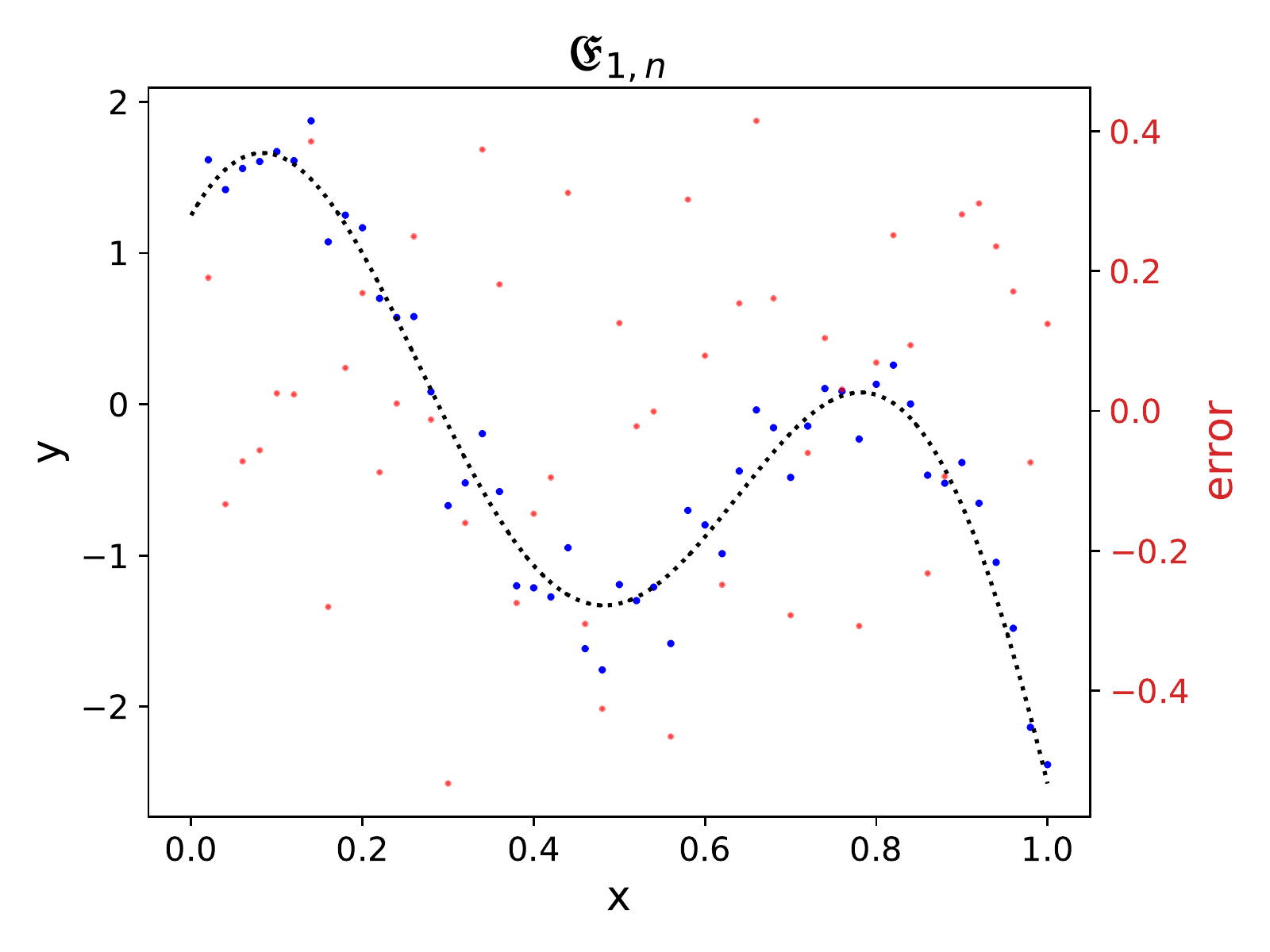}
\end{subfigure}
     \begin{subfigure}[b]{0.49\textwidth}
         \centering
         \includegraphics[width=\textwidth]{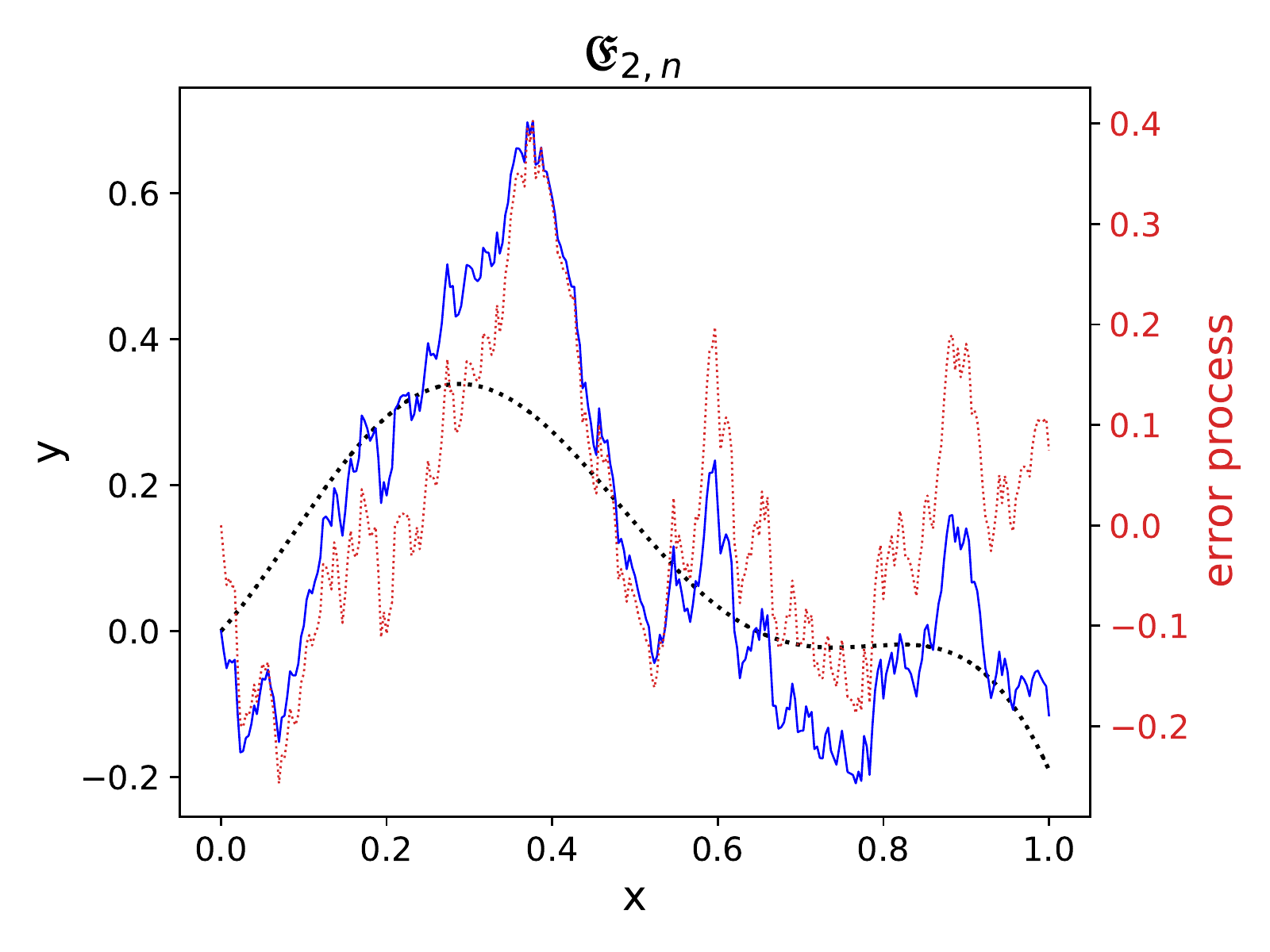}
\end{subfigure}
    \caption{Exemplary realizations of the observations in the two experiments $\Expm_{1,n}$ (left) and $\Expm_{2,n}$ (right) for $n=50$ where the Gauss-Markov process is given through the functions $u(x)=\exp(x)$ and $q(x) = \exp(2x) - 1$. In both figures the actual observations are plotted in \textcolor{blue}{blue}. The black dotted line represents the true function $f$ in the left plot, and its anti-derivative $F_f(\cdotb)=\int_0^{\cdotb} f(s)\dd s$ in the right plot.
    The error terms, that is, the increments of $\Xi$ for $\Expm_{1,n}$ and the process $\Xi$ itself for $\Expm_{2,n}$, are plotted in \textcolor{red}{red} (note that there is a secondary y-axis in \textcolor{red}{red} for the error terms).}
    \label{fig:plots}
\end{figure}

\noindent The second {experiment with continuous observations} is given through
\begin{equation*}
	\Expm_{2,n}(\Theta) = (\Cc([0,1], \R), \Bs(\Cc([0,1], \R)), (\Pb_{2,n}^f)_{f \in \Theta})~,
\end{equation*}
where $\Pb_{2,n}^f$ denotes the distribution of the continuous time process $(Y_t)_{t \in [0,1]}$ defined by 
\begin{equation}\label{eq:Y_t}
	Y_t = \int_0^t f(s) \dd s + \frac{1}{\sqrt n} \Xi_t, \quad t \in [0,1].
\end{equation}
Typical realizations from the two experiments are visualized in Figure~\ref{fig:plots} for the case of an Ornstein-Uhlenbeck process conditioned to start in zero.

The following  result provides sufficient conditions for asymptotic equivalence of the experiments $\Expm_{1,n}$ and $\Expm_{2,n}$ in terms of the class $\Theta$ and the Gauss-Markov process $\Xi$. The proof   is given in Appendix~\ref{app:general}.

\begin{theorem}\label{thm:general:triangular}
	Let $\Xi$ be a Gauss-Markov process on the interval $[0,1]$ with $\Xi_0 = 0$ such that Assumption~\ref{ass:GM} holds with $K_\Xi(s,t) = u(s)v(t)$ for $0 \leq s \leq t \leq 1$, and $v(1) \neq 0$.
	Set $F_{f}(\boldsymbol\cdotb) = \int_0^{\boldsymbol\cdotb} f(s)\dd s$ and assume that the conditions 
\begin{enumerate}[(i)]
		\item\label{triC1} $\frac{1}{n} \sup_{f \in \Theta} \sum_{i=1}^{n} \frac{(f(t_{i,n})  - n \int_{t_{i-1,n}}^{t_{i,n}} f(s)\dd s)^2}{v^2(t_{i,n}) (q(t_{i,n}) - q(t_{i-1,n}))} \to 0$\qquad and
		\bigskip
		\item\label{triC2} $\sqrt{n} \sup_{f \in \FctClass} \inf_{\alphab_n \in \R^n} \lVert F_{f}(\cdotb) - \sum_{i=1}^n \alpha_i K_\Xi\left( \cdotb, t_{i,n} \right) \rVert_{\RKHS(\Xi)} \to 0$
	\end{enumerate}
\medskip
are satisfied, 
	where $t_{i,n} = i/n$ and $\lVert \, \cdotb \, \rVert_{\RKHS(\Xi)}$ is the norm in the RKHS $\RKHS(\Xi)$.
	Then, the two sequences of experiments $(\Expm_{1,n}(\Theta) )_{n \in \N}$ and $(\Expm_{2,n}(\Theta) )_{n \in \N}$ are asymptotically equivalent, i.e.
	\begin{equation*}
		\Expm_{1,n}(\Theta) \approx \Expm_{2,n}(\Theta).
	\end{equation*}
\end{theorem}

\section{Specific function classes and a counterexample}\label{s:sobolev}

In this section, we first apply the abstract Theorem~\ref{thm:general:triangular} in order to establish asymptotic equivalence over Sobolev ellipsoids (Section~\ref{subs:sobolev}) and Hölder balls (Section~\ref{subs:hoelder}) of sufficiently smooth functions.
In Section~\ref{subs:sketch} we discuss an alternative proof strategy, which is based on some transformations of the different experiments and requires stronger assumptions as used in Section \ref{subs:sobolev} and \ref{subs:hoelder}.
Finally, in Section~\ref{s:BB} we give a counterexample which shows that asymptotic equivalence cannot hold (even for function classes containing very smooth functions) in general if only Assumption~\ref{ass:GM} is required.

\subsection{Asymptotic equivalence for Sobolev ellipsoids}\label{subs:sobolev}
Consider the (complex) trigonometric basis $\{ \eb_k \}_{k \in \Z}$ of $L^2([0,1])$ given by
\begin{equation*}
	\eb_k (\cdotb) = \exp(-2\pi \ii k \cdotb).
\end{equation*}
Any function $f$ in $L^2([0,1])$ can be represented as a convergent series (in $L^2$-sense) of the form
\begin{equation*}
	f(\cdotb) = \sum_{k=-\infty}^\infty \theta_k \eb_k(\cdotb)
\end{equation*}
where $\theta_k = \langle f,\eb_k \rangle_{L^2} = \int_0^1 f(s)\exp(2\pi \ii k s)\dd s$ are the uniquely defined Fourier coefficients of the function $f$.
As potential function classes over which we would like to establish asymptotic equivalence we consider Sobolev ellipsoids
\begin{equation} \label{sob} 
	\Theta(\beta,L) = \left\lbrace f(\cdotb) = \sum_{k=-\infty}^\infty \theta_k \eb_k(\cdotb) : \theta_k  = \overline{\theta_{-k}} \text{ and } \sum_{k=-\infty}^\infty (1 + \lvert k \rvert)^{2\beta} \lvert \theta_k \rvert^2 \leq L^2 \right\rbrace. 
\end{equation}
Note that the condition that $\theta_k = \overline{\theta_{-k}}$ is used only to ensure that the functions  in the set 
$\Theta(\beta,L)$ are real-valued.
The following result establishes asymptotic equivalence of the experiments $\Expm_{1,n}(\Theta(\beta,L))$ and $\Expm_{2,n}(\Theta(\beta,L))$ considered in Section~\ref{s:cond} for $\beta > {1}/{2}$ under some assumptions on the functions $v$ and $q$
characterizing the Gauss-Markov process $(\Xi_t)_{t\in [0,1]}$.

\begin{theorem}\label{thm:sobolev}
	Let $\Xi$ be a Gauss-Markov process on the interval  $[0,1]$ with $\Xi_0 = 0$ such that Assumption~\ref{ass:GM} holds with $K_\Xi(s,t) = u(s)v(t)$ for $0 \leq s \leq t \leq 1$. If 
 $\beta > 1/2$ and the assumptions 
	\begin{enumerate}[a)]
		\item\label{it:sob:a} $q \in C^1$ and $0 < \inf_{t \in [0,1]} q'(t) \leq \sup_{t \in [0,1]} q'(t) \leq C < \infty$,
		\item\label{it:sob:b} $\inf_{t \in [0,1]} v(t) > 0$, and
		\item\label{it:sob:c} $v'$ and $q'$ are Hölder continuous with index $\gamma>1/2$
	\end{enumerate}
are satisfied, then 
	\begin{equation*}
		\Expm_{1,n}(\Theta(\beta,L)) \approx \Expm_{2,n}(\Theta(\beta,L)).
	\end{equation*}
\end{theorem}
The proof of Theorem~\ref{thm:sobolev} is complicated and deferred to Appendix~\ref{app:sobolev}.
We proceed with some remarks that put the statement  in the context of well-known results for the  Brownian motion.

\begin{remark}
	It is known that for the smoothness index $\beta = 1/2$ asymptotic equivalence between nonparametric regression model \eqref{eq:BL:regr} and the white noise with drift experiment  \eqref{eq:BL:GWN} (i.e., our experiment $\Expm_{2,n}$ for the special case of Brownian motion) does not hold. See Remark~4.6 in \cite{brown1996asymptotic} for further details.
	In Section~\ref{s:BB} below  we consider the example of a Brownian bridge as an error process and  prove that asymptotic equivalence for arbitrary values of $\beta$ cannot hold without additional assumptions on the Gauss-Markov process (like, for instance, assumptions on the functions $v$ and $q$ as made in Theorem~\ref{thm:sobolev}).
\end{remark}

\begin{remark}
	In the proof of Theorem~\ref{thm:sobolev} it is shown that under our assumptions the Le Cam distance between the experiments $\Expm_{1,n}(\Theta(\beta,L))$ and $\Expm_{2,n}(\Theta(\beta,L))$ converges to zero with a rate that cannot be faster than
	\begin{equation*}
	    \Oc(\max \{ n^{-1/2},n^{-\beta + 1/2} \})
	\end{equation*}
	(indeed, the rate might even be slower due to the terms incorporating the H\"older index $\gamma$ from Condition \ref{it:sob:c} in the statement of the theorem).
	In the case where the Gauss-Markov is standard Brownian motion the rate $n^{-\beta + 1/2}$ can be established which is faster for $\beta > 1$; see Rohde, \cite{rohde2004asymptotic} for further details.
	Since we have restricted ourselves to proving asymptotic equivalence, we have not further pursued this issue here.
\end{remark}

\begin{remark}[Sequence space representations]
	Note that for the special case of Brownian motion the relation between the functions $f=F^\prime$ and $g$ in  in the representation 
	\eqref{eq:RKHS:explicit} of $F \in \RKHS(\Xi) $  reduces to  the equality 
	$g=f$.
	Hence, by representing $f$ in terms of the  coefficients of a series expansion with respect to an orthonormal basis (for instance, the complex trigonometric basis given through $\{ \eb_k \}_{k \in \Z}$) the same representation is achieved for $g$.
	Via the isomorphism $\psi$  defined in \eqref{psi} (more precisely, its inverse $\psi^{-1}$) this yields a representation of $F_f$ in terms of an orthonormal basis in the RKHS leading to a Gaussian sequence space model of the form
	\begin{equation*}
		y_k = \theta_k + \xi_k, \qquad k \in \Z,
	\end{equation*}
	where $\xi_k$ i.i.d.\,$\sim \Nc(0,1)$ and the sequence $\theta = (\theta_k)_{k \in \Z}$ satisfies the same smoothness conditions as in the definition of the Sobolev ellipsoid  $\Theta(\beta,L)$.
	In the general case of an arbitrary Gauss-Markov process the relation \eqref{eq:g:F} does not 'transport' an orthonormal basis for the representation of $f$ to an orthonormal basis for $g$.

	\bigskip

	An alternative sequence space representation (for general Gaussian processes) is obtained from \eqref{eq:Y_t} using the Karhunen-Loeve representation of the process $\Xi$.
	Note that the kernel $K_\Xi$ gives rise to an integral operator defined through
    \begin{equation*}(K_\Xi f)(s) = \int_0^1 K_\Xi(s,t)f(t)\dd t
    \end{equation*}
    for $f \in L^2([0,1])$ (by slight abuse of notation, we use the same letter for both the kernel and the operator associated with it).
    Provided that (the kernel) $K_\Xi$ is continuous, Mercer's theorem applies and (the operator) $K_\Xi$ has a countable system of eigenfunctions $\phi_k$ with corresponding eigenvalues $\lambda_k \geq 0$.
    By Proposition~3.11 in Neveu \cite{neveu1968processus} an orthonormal basis of the RKHS $\RKHS(\Xi)$ is given by the functions
    $g_k = \sqrt{\lambda_k} \phi_k$, $k=1,2,\ldots$ with $\lambda_k > 0$.
    In terms of this basis, \eqref{eq:Y_t} can be rewritten as
    \begin{equation}\label{eq:seq:model}
        y_k = \theta_k + \frac{1}{\sqrt{n}} \eta_k, \quad k = 1,2,\ldots
    \end{equation}
    where $\theta_1, \theta_2 , \ldots $ are the Fourier coefficients of the function $F$ with respect to the basis $g_1, g_2, \ldots $
 and $\eta_1, \eta_2, \ldots  $ are independent and standard normal distributed  (note that the sequence model here is stated in terms of the coefficients of $F_f$ instead of those of $f$ itself).
    In the special case of Brownian motion, one has (see Beder, \cite{beder1987sieve}, p.~66)
    \begin{equation*}
    g_k(t) = \frac{\sqrt{2}}{(k-\frac{1}{2})\pi} \sin \left( (k-\frac{1}{2})\pi t \right) = \sqrt{\lambda_k} \phi_k, \quad \lambda_k = \frac{1}{(k-\frac{1}{2})^2 \pi^2}.
    \end{equation*}
    In this case, the derivatives $\psi_k = g_k'$ of the basis functions $g_k$ in the RKHS form an orthonormal basis of $L^2([0,1])$.
    Hence, the sequence space model \eqref{eq:seq:model} can be interpreted both as a perturbation of the coefficients of $F_f$ or of $f$ (note that $f \in L^2([0,1]) $ if and only if $  F_f \in \RKHS(\Xi)$).

\end{remark}

\subsection{Asymptotic equivalence for H\"older classes}\label{subs:hoelder}

For a smoothness index $0 < \alpha \leq 1$ and  constants $0 < L < \infty$, $0 < M \le \infty$, we introduce the H\"older class
\begin{equation*}
    \Fc(\alpha, L, M) = \left\{ f \colon [0,1] \to \R : \sup_{0 \leq x < y \leq 1} \frac{\lvert f(x) - f(y) \rvert}{\lvert x - y \rvert^\alpha} \leq L \text{ and } \lVert f \rVert_\infty \leq M \right\}~,
\end{equation*}
Where the case $M=\infty$ means that the assumption  $ \lVert f \rVert_\infty \leq M $ is omitted in the definition of 
$ \Fc(\alpha, L, M) $. 
Under the same technical assumptions on the functions $v$ and $q$ as  stated in Theorem~\ref{thm:sobolev}  we obtain asymptotic equivalence for functions from H\"older classes with  smoothness  index $\alpha > 1/2$.
\begin{theorem}\label{thm:hoelder}
	Let $\Xi$ be a Gauss-Markov process on the interval  $[0,1]$ with $\Xi_0 = 0$ such that Assumption~\ref{ass:GM} holds with $K_\Xi(s,t) = u(s)v(t)$ for $0 \leq s \leq t \leq 1$. Suppose further that 
  $\alpha > 1/2$ and  that the assumptions 
	\begin{enumerate}[a)]
		\item $q \in C^1$ and $0 < \inf_{t \in [0,1]} q'(t) \leq \sup_{t \in [0,1]} q'(t) \leq C < \infty$,
		\item $\inf_{t \in [0,1]} v(t) > 0$, and
		\item $v'$ and $q'$ are Hölder continuous with index $\gamma >1/2$
	\end{enumerate}
	are satisfied.
	Then, if $ M < \infty$, we have
	\begin{equation} \label{hd1} 
		\Expm_{1,n}(\Fc(\alpha,L,M)) \approx \Expm_{2,n}(\Fc(\alpha,L,M)).
	\end{equation}
	Moreover, if the function $v\cdot q'$ is constant on the interval $[0,1]$, the  asymptotic equivalence \eqref{hd1} also holds in the case $M=\infty$.
\end{theorem}
The proof of Theorem~\ref{thm:hoelder} is given in Appendix~\ref{app:hoelder}.
Note that in the Gaussian white noise model we have $v\cdot q' \equiv 1$, and in this case Theorem~\ref{thm:hoelder} reproduces the result given in \cite{brown1996asymptotic}, at the top of page~2390.

\subsection{Sketch of a different proof approach}\label{subs:sketch}

The proofs of our results concerning Sobolev and Hölder classes in Sections~\ref{subs:sobolev} and \ref{subs:hoelder} both relied on the general conditions stated in Theorem~\ref{thm:general:triangular}.
We will now sketch a different proof that is not based on the RKHS approach 
used in this paper.
It requires a slightly stronger assumption on the function $q$ and exploits results from the seminal paper \cite{brown1996asymptotic}, which are applied to suitable transformations of the experiments of interest.
Although checking all the details of this alternative approach might be a tedious affair, we will state the
main steps here which yields a look on
the asymptotic equivalence from a different angle.

The arguments used in the sequel make use of several transformations of the experiments.
To be precise, starting  with the discrete experiment $\Expm_{1,n}$, we note that observing the $Y_{i,n}$ in Equation~\eqref{eq:Yin} is equivalent to observing
\begin{equation*}
    Y^\star_{j,n} = \frac{\sum_{i=1}^j Y_{i,n}}{v(t_{j,n})} = \frac{\sum_{i=1}^j f(t_{i,n})}{v(t_{j,n})} + \sqrt{n} W_{q(t_{j,n})}~,~~j = 1,\ldots,n. 
\end{equation*}
By taking increments these observations are in turn equivalent to the experiment $\Expm_{1,n}^\star$ with observations
\begin{equation}\label{eq:Yjn:star}
   Y^\star_{j,n} - Y^\star_{j-1,n} = \mu_{j,n} + \sigma_{j,n} \epsilon_{j,n} ~,~~
   j=1,\ldots,n,
\end{equation}
where $Y^\star_{0,n} = 0$, 
the $\epsilon_{1,n} ,\ldots,\epsilon_{n,n}$ are i.i.d.\,$\sim \Nc(0,1)$,
\begin{align*}
    \mu_{j,n} = \frac{\sum_{i=1}^j f(t_{i,n})}{v(t_{j,n})} - \frac{\sum_{i=1}^{j-1} f(t_{i,n})}{v(t_{j-1,n})},
\end{align*}
and
\begin{align*}
    \sigma_{j,n}^2 = n \left[ q(t_{j,n}) - q(t_{j-1,n}) \right].
\end{align*}
We now transform the continuous experiment $\Expm_{2,n}$ to observations with a similar structure as in \eqref{eq:Yjn:star} but different mean and variance components.
The observation of the path $Y_t$ in \eqref{eq:Y_t} is, assuming that the function $v$ does not vanish, equivalent to observing the path
\begin{equation*}
    Z_t^\star = \frac{Y(t)}{v(t)} = \int_0^t \mu(s) \dd s + \frac{1}{\sqrt n} W_{q(t)}, \qquad t \in [0,1].
\end{equation*}
Here, the function $\mu\colon [0,1]\to \R$ is given through
\begin{equation*}
    \mu(s) = \frac{f(s)v(s) - v'(s) \int_0^s f(\tau)\dd \tau}{v^2(s)}.
\end{equation*}
Using time-change for martingales (see \cite{karatzas1991brownian}, Chapter~3, Theorem~4.6) the Gaussian process $Z_t^\star$ can be characterized by the SDE
\begin{equation}\label{eq:alter:2:cont}
    \dd Z^\star(t) = \mu(t) \dd(t) + \frac{1}{\sqrt{n}} \sqrt{q'(t)} \dd W_t. 
\end{equation}
Now, under regularity assumptions on $\mu$ and $q$  Theorem 4.1 in \cite{brown1996asymptotic} can be applied to show asymptotic equivalence of the experiments where one either observes the path $Z_t^\star$ or the discrete data
\begin{equation}\label{eq:alter:2:disc}
    Z_{j,n}^\star = \mu(j/(n+1)) + \sqrt{q'(j/(n+1))} \epsilon_{j,n}~,~~j=1,\ldots,n,
\end{equation}
and  $\epsilon_{1,n}, \ldots, \epsilon_{n,n}$ are i.i.d.\,$\sim \Nc(0,1)$.
Let us denote the latter  experiment by 
$\Expm_{2,n}^\star$.
Then,  by the arguments  sketched  above the target equivalence $\Expm_{1,n} \approx \Expm_{2,n}$ follows from the equivalence $\Expm_{1,n}^\star \approx \Expm_{2,n}^\star$.
In order to prove this statement, however, it is sufficient to show that
\begin{equation*}
    n \cdot \sup_j\left( \mu (j/(n+1)) - \mu_{j,n} \right)^2 + n \sup_j \left( q'(j/(n+1)) - \sigma_{j,n}^2 \right)^2 \to 0,
\end{equation*}
which can be established under mild regularity assumptions.

We emphasize that the assumptions for this 
alternative proof strategy and the ones imposed by us in Theorems \ref{thm:sobolev} and \ref{thm:hoelder}, respectively, do not coincide entirely.
In particular, Theorems~\ref{thm:sobolev} and \ref{thm:hoelder} in the paper assume the existence of a
 Hölder continuous
first derivative $q'$. On the other hand  
Theorem~4.1 in \cite{brown1996asymptotic}, which is   used to establish  the asymptotic equivalence of the models~\eqref{eq:alter:2:cont} and \eqref{eq:alter:2:disc} in this section, requires 
the existence of a second derivative of $q$ (see Condition~(4.1) in this reference). 
From this point of view, the RKHS approach uses  
(slightly) weaker assumptions.

\subsection{Non-equivalence for Brownian bridge}\label{s:BB}

We conclude this section by showing that the asymptotic equivalence result established in Theorem~\ref{thm:sobolev} cannot be transferred to arbitrary Gauss-Markov processes.
To be more precise, Theorem~\ref{thm:bridge} below shows that for the special case of the Brownian bridge the deficiency $\delta(\Expm_{1,n},\Expm_{2,n})$ is bounded away from zero, and as a consequence the two experiments are not asymptotically equivalent over Sobolev spaces $\Theta(\beta,L)$ for arbitrary smoothness $\beta$ and $L > 0$. Recall that the covariance kernel of a Brownian Bridge $(B_t)_{t \in [0,1]}$ satisfies the representation
\begin{equation*}
    \Cov(B(s),B(t)) = u(s) v(t), \qquad 0 \leq s \leq t \leq 1,
\end{equation*}
with $u(s)=s$, $v(t)=1-t$.
Therefore $v(1)=0$, and in addition Conditions \ref{it:sob:a} and \ref{it:sob:b} in Theorem \ref{thm:sobolev} are not satisfied.
\begin{theorem}\label{thm:bridge}Consider the two sequences of experiments $(\Expm_{1,n})_{n\in \N} $ and $(\Expm_{2,n})_{n\in \N} $  defined  in \eqref{eq:Yin} and \eqref{eq:Y_t}, respectively, where the error process $(\Xi_t)_{t\in [0,1]}$ is a Brownian Bridge.
	Then, for an arbitrary smoothness index $\beta > 0$, any $L > 0$, and all $n \in \N$, 
	\begin{equation*}
		\delta(\Expm_{1,n}(\Theta(\beta,L)),\Expm_{2,n}(\Theta(\beta,L))) \geq 1/4,
	\end{equation*}
	where $\Theta(\beta,L)$ is the Sobolev ellipsoid introduced in	\eqref{sob}.
	In particular, the experiments $(\Expm_{1,n})_{n\in \N} $ and $(\Expm_{2,n})_{n\in \N} $  are \emph{not} asymptotically equivalent.
\end{theorem} 
\newpage

\appendix

\section{Proof of Theorem~\ref{thm:general:triangular}}\label{app:general}

The proof consists in the consideration of two intermediate experiments, given through Equations \eqref{eq:Yinp} and \eqref{eq:triKriging} below, that lie between $\Expm_{1,n}$ and $\Expm_{2,n}$.

\medskip

\noindent\underline{First step}: We first show that under Condition~\ref{triC1} in Theorem~\ref{thm:general:triangular}  the model given by Equation~\eqref{eq:Yin} is asymptotically equivalent to observing
\begin{align}\label{eq:Yinp}
	Y_{i,n}' = n \int_{t_{i-1,n}}^{t_{i,n}} f(s)\dd s + \sqrt{n} \xi_{i,n}, \quad i=1,\ldots,n,
\end{align}
where $\xi_{i,n}= \Xi_{t_i,n} -\Xi_{t_{i-1,n}}$ are the increments of the process $(\Xi_t)_{t\in [0,1]}$.
Under Assumption~\ref{ass:GM}, we can take advantage of the representation \eqref{eq:repr:GM} and write
\begin{align*}
	\xi_{i,n} 
&=v(t_{i,n}) W_{q(t_{i,n})} - v(t_{i-1,n})W_{q(t_{i-1,n})},
\end{align*}
and the experiment \eqref{eq:Yin} can be written as
\begin{equation*}Y_{i,n} = f ( t_{i,n}) + \sqrt{n} \left[ v(t_{i,n}) W_{q(t_{i,n})} - v(t_{i-1,n})W_{q(t_{i-1,n})} \right].
\end{equation*}
Adding and subtracting $v(t_{i,n}) W_{q(t_{i-1,n})}$, we get
\begin{align*}
	Y_{i,n} & = f(t_{i,n}) + \sqrt{n} (v(t_{i,n}) - v(t_{i-1,n})) W_{q(t_{i-1,n})}  \\ 
	&  ~~~~~~~~~~~~~ 	+ \sqrt{n} v(t_{i,n}) (W_{q(t_{i,n})} - W_{q(t_{i-1,n}))}.
\end{align*}
Similarly, \eqref{eq:Yinp} can be written as 
\begin{align*}
	Y_{i,n}' &= n \int_{t_{i-1,n}}^{t_{i,n}} f(s)\dd s + \sqrt{n} (v(t_{i,n}) - v(t_{i-1,n})) W_{q(t_{i-1,n})}  \\
	&  ~~~~~~ ~~~~~~ ~~~~~~ ~~~~~~ + \sqrt{n} v(t_{i,n}) (W_{q(t_{i,n})} - W_{q(t_{i-1,n})}).
\end{align*}
For the sake of a transparent notation let $\Pb^{\Yb_n} = \Pb_{1,n}^f $ denote the distribution of the vector $\Yb_n=(Y_{1,n},\ldots,Y_{n,n})$, where we do not reflect the dependence on $f$ in the notation. Similarly, let 
$\Pb^{\Yb_n'}$ denote the distribution of the vector $\Yb_n' = (Y'_{1,n},\ldots,Y'_{n,n})$.
Note  that the squared total variation distance can be bounded by the Kullback-Leibler divergence, and therefore we will  derive a
 bound for the Kullback-Leibler distance between the distributions $\Pb^{\Yb_n}$ and $\Pb^{\Yb_n'}$   by suitable conditioning.
Denote with $\Fs_{i,n}$ the $\sigma$-algebra generated by $\left\lbrace W_{q({t})}, t \leq t_{i,n} \right\rbrace$.
We have
\begin{align*}
	\KL(\Pb^{\Yb_n}, \Pb^{\Yb_n'}) &= \Eb [  \KL( \Pb^{\Yb_n}, \Pb^{\Yb_n'} | \Fs_{n-1,n}) ] \\
	&= \Eb [ \KL (\Pb^{\Yb_{1:n-1,n}}, \Pb^{\Yb_{1:n-1,n}'} | \Fs_{n-1,n}) +  \KL (\Pb^{Y_{n,n}}, \Pb^{Y_{n,n}'} | \Fs_{n-1,n})]\\
	&= \Eb [ \KL (\Pb^{\Yb_{1:n-1,n}}, \Pb^{\Yb_{1:n-1,n}'} | \Fs_{n-1,n}) ] + \Eb [ \KL (\Pb^{Y_{n,n}}, \Pb^{Y_{n,n}'} | \Fs_{n-1,n})]\\
	&=  \KL(\Pb^{\Yb_{1:n-1,n}}, \Pb^{\Yb_{1:n-1,n}'})  + \Eb [ \KL (\Pb^{Y_{n,n}}, \Pb^{Y_{n,n}'} | \Fs_{n-1,n})] ,
\end{align*}
where we use the notation $\Yb_{1:n-1,n}=(Y_{1,n},\ldots,Y_{n-1,n})$  and  $\Yb_{1:n-1,n}' = (Y'_{1,n},\ldots,Y'_{n-1,n})$. 
Repeating the same argument, one obtains
\begin{equation*}
	\KL(\Pb^{\Yb_n}, \Pb^{\Yb_n'}) = \sum_{i=1}^n \Eb [ \KL (\Pb^{Y_{i,n}}, \Pb^{Y_{i,n}'} | \Fs_{i-1,n})] ].
\end{equation*}
In order to study the terms  $\Eb [ \KL (\Pb^{Y_{i,n}}, \Pb^{Y_{i,n}'} | \Fs_{i-1,n})]$, note that 
\begin{align*}
	Y_{i,n} | \Fs_{i-1,n} \sim \Nc ( \mu_i , \sigma_i^2 ) \qquad \text{and} \qquad Y'_{i,n} | \Fs_{i-1,n} \sim \Nc ( \mu_i^\prime, \sigma_i^2 ),
\end{align*}
where
\begin{align*}
    &\mu_i = f(t_{i,n}) + \sqrt{n} (v(t_{i,n}) - v(t_{i-1,n})) W_{q(t_{i-1,n})},\\
    &\mu_i^\prime = n \int_{t_{i-1,n}}^{t_{i,n}} f(s)\dd s + \sqrt{n} (v(t_{i,n}) - v(t_{i-1,n})) W_{q(t_{i-1,n})},\quad \text{and}\\
    &\sigma_i^2 = n v^2(t_{i,n}) (q(t_{i,n}) - q(t_{i-1,n})).
\end{align*}
Here and in the following, $\Nc(\mu , \sigma^2)$ denotes a normal distribution with mean $\mu$ and variance $\sigma^2$.
Using the fact that the Kullback-Leibler divergence between two normal distributions with common variance is given by
\begin{equation*}
	\KL(\Nc(\mu_1, \sigma^2), \Nc(\mu_2, \sigma^2)) = \frac{(\mu_1 - \mu_2)^2}{2\sigma^2},
\end{equation*}
we have
\begin{equation*}
	\KL (\Pb^{Y_{i,n}}, \Pb^{Y_{i,n}'} | \Fc_{i-1,n}) = \frac{(f(t_{i,n})  - n \int_{t_{i-1,n}}^{t_{i,n}} f(s)\dd s)^2}{2 nv^2(t_{i,n}) (q(t_{i,n}) - q(t_{i-1,n}))}.
\end{equation*}
This yields
\begin{equation*}
\KL(\Pb^{\Yb_n}, \Pb^{\Yb_n'}) = \frac{1}{2n} \sum_{i=1}^{n} \frac{(f(t_{i,n})  - n \int_{t_{i-1,n}}^{t_{i,n}} f(s)\dd s)^2}{v^2(t_{i,n}) (q(t_{i,n}) - q(t_{i-1,n}))},
\end{equation*}
and $\KL(\Pb^{\Yb_n}, \Pb^{\Yb_n'}) \to 0$ holds if and only if Condition~\ref{triC1} holds.
Consequently, the experiments \eqref{eq:Yin} and \eqref{eq:Yinp} are asymptotically equivalent in this case. 
\medskip

\noindent\underline{Second step}:
Let $I(t | \yb_n)$ denote the Kriging interpolator which is defined as
\begin{equation}\label{eq:triKriging}
I(t| \yb_n) = \left(  K_\Xi ( t,t_{1,n}) , K_\Xi \left( t,t_{2,n}\right) , \ldots, K_\Xi ( t,t_{n,n})  \right) \Cov(\Xib_n)^{-1} \yb_n^\top
\end{equation}
for $\yb_n = (y_{1,n},\ldots,y_{n,n})$ and $\Xib_n = (\Xi_{t_{1,n}},\ldots,\Xi_{t_{n,n}})$ (the additional condition that $v(1) \neq 0$ guarantees the invertibility of $\Xib_n$; see Lemma~A.1 in \cite{dette2016optimal} for an explicit formula for the entries of the inverse matrix).
By definition 
the Kriging predictor is linear in the argument $\yb_n$, and a simple argument shows the interpolation property 
\begin{equation}
    \label{hd2} 
    I(t_{j,n}|\yb_n) = y_{j,n} ~\text{  for } j =1,\ldots, n.
\end{equation}
The second step now consists in proving (exact, that is, non-asymptotic) equivalence of the experiment defined by the discrete observations \eqref{eq:Yinp} and the experiment defined by the continuous  path
\begin{equation}\label{eq:path:Kriging}
\Ytilde_t = I(t|\Fb_{f,n}) + n^{-1/2} \Xi_t,
\end{equation}
where $\Fb_{f,n}=(F_{f,n}(t_{1,n}),\ldots,F_{f,n}(t_{n,n}))$.
Defining the partial sums  $S'_{k,n} = \sum_{j=1}^{k} Y'_{j,n}$ and recalling the notation
$\xi_{k,n}= \Xi_{t_{k,n}} -\Xi_{t_{k-1,n}}  $ for the increments of the the process $(\Xi_t)_{t\in [0,1]}$, we have
\begin{align} \label{hd3}
	S'_{k,n} = \sum_{j=1}^{k} Y'_{j,n} &= n \int_0^{t_{k,n}} f(s) \dd s + \sqrt{n} \sum_{j=1}^{k} \gn_{j,n} 
	= n F_f\left( t_{k,n} \right) + \sqrt{n} \Xi_{t_{k,n}}
\end{align}
where we used the interpolating property \eqref{hd2} and the the definition \eqref{eq:path:Kriging}.
Let $(\Xi_t^\prime)_{t \in [0,1]}$ be an independent copy of $(\Xi_t)_{t \in [0,1]}$, and set $R_t = \Xi_t^\prime - I(t| \Xib_n^\prime)$ with $\Xib_n^\prime = (\Xi^\prime_{t_{1,n}},\ldots,\Xi^\prime_{t_{n,n}})$.
Then, the process
\begin{equation*}
	I(t|\Xib_n) + R_t, \quad t \in [0,1],
\end{equation*}
follows the same law as $(\Xi_t)_{t \in [0,1]}$ and  $(\Xi_t^\prime)_{t \in [0,1]}$, which can can be checked by a comparison of the covariance structure (indeed, this kind of construction is valid for any centered Gaussian process).
Then, observing the definition \eqref{eq:path:Kriging}, we have 
\begin{align*}
	\Ytilde_t &= 
	I(t|\Fb_{f,n}) + n^{-1/2} \Xi_t \\
	&\stackrel{\Lc}{=}  I(t|\Fb_{f,n}) + n^{-1/2}I(t|\Xib_n) + n^{-1/2}R_t \\
	&= 	n^{- 1} (I(t|\Sb'_n) + \sqrt n R_t),
\end{align*}
where we used the notation $\Sb'_n = (S'_{1.n}, \ldots , S'_{n.n})$,  equation \eqref{hd3} and the linearity of the Kriging estimator.
Therefore, the process $(\Ytilde_t)_{t \in [0,1]}$ can be constructed  from the vector $\Yb'_{n}$.
On the other hand,  the observations $Y'_{1,n},\ldots,Y'_{n,n}$ can be recovered from the trajectory $(\Ytilde_t)_{t \in [0,1]}$ since for $t=t_{k,n}$   the interpolation property \eqref{hd2} yields 
\begin{align*}
	n \Ytilde_{t_{k,n}} &= n F_{f,n}\left( t_{k,n} \right) + \sqrt n \Xi_{t_{k,n}}\\
	&= n \int_0^{t_{k,n}} f(s) \dd s + \sqrt n(\Xi_{t_{k,n}} - \Xi_{0})\\
	&=  n \int_0^{t_{k,n}} f(s) \dd s + \sqrt n \sum_{j=1}^k \xi_{j,n},
\end{align*}
and one obtains $Y'_{i,n}$ as $Y'_{i,n} =n\Ytilde_{t_i} - n\Ytilde_{t_{i-1}}$.
Hence, the process  $(\Ytilde_t)_{t \in [0,1]}$ and  the vector $\Yb'_n$ contain the same information 
and the experiments \eqref{eq:Yinp} and  \eqref{eq:path:Kriging} are equivalent.
\medskip

\noindent\underline{Third step}:
It remains to show that the experiment  $\widetilde{\Expm}_{2,n}$ defined by the path  in \eqref{eq:path:Kriging} is asymptotically equivalent to 
the experiment $\Expm_{2,n}$ defined by  
\begin{equation*}
	Y_t = \int_0^t f(s)\dd s + n^{-1/2} \Xi_t.
\end{equation*}
For this purpose we denote by $ \Pb^{(Y_t)_{t \in [0,1]}}$ and 
$ \Pb^{(\Ytilde_t)_{t \in [0,1]} }$ the distributions of
the processes  $(Y_t)_{t \in [0,1]}$  and $(\Ytilde_t)_{t \in [0,1]}$, respectively, where the dependence 
on the parameter $f$ is again suppressed.
First, note that the representation of the Kriging estimator shows that the function $t \to I(  t | \Fb_{f,n})$ belongs to the RKHS $\RKHS(\Xi)$ associated with the covariance kernel $K_\Xi$.
Second, Condition~\ref{triC2} in the statement of Theorem~\ref{thm:general:triangular} yields that the same holds true for the function $F_{f,n}(\boldsymbol \cdot)$.
Using the fact that the squared total variation distance can be bounded by the Kullback-Leibler divergence, we obtain
\begin{align*}
	\Delta^2(\Expm_{2,n}, \widetilde \Expm_{2,n}) &\leq \sup_{f \in \Theta} \KL( \Pb^{(Y_t)_{t \in [0,1]}}, \Pb^{(\Ytilde_t)_{t \in [0,1]}})\\
	&= \frac{n}{2} \sup_{f \in \FctClass} \lVert F_{f}(\boldsymbol\cdot) - I(\boldsymbol\cdot \mid \Fb_{f,n}) \rVert_{\RKHS(\Xi)}^2\\
	&= \frac{n}{2} \sup_{f \in \FctClass} \inf_{\alphab_n \in \R^n} \left\lVert F_f(\boldsymbol\cdotb) - \sum_{j=1}^n \alpha_j K_\Xi \left( \boldsymbol\cdot, t_{j,n} \right)  \right\rVert_{\RKHS(\Xi)}^2 \to 0
\end{align*}
(the first equality follows from Lemma~2 in \cite{schmidt-hieber2014asymptotic}, the second from Theorem~13.1 in \cite{wendland2005scattered}), which completes the proof of Theorem~\ref{thm:general:triangular}. 
\section{Proofs of the results in Section~\ref{s:sobolev}}

\subsection{Proof of Theorem~\ref{thm:sobolev}}\label{app:sobolev}

The proof consists in checking the two conditions \ref{triC1} and \ref{triC2} in Theorem~\ref{thm:general:triangular}.

\medskip

\noindent\underline{Verification of condition \ref{triC1}}:
We have to show that the expression
\begin{equation*}
	\frac{1}{n} \sup_{f \in \Theta(\beta,L)} \sum_{i=1}^{n} \frac{(f(t_{i,n})  - n \int_{t_{i-1,n}}^{t_{i,n}} f(s)\dd s)^2}{v^2(t_{i,n}) (q(t_{i,n}) - q(t_{i-1,n}))}
\end{equation*}
converges to zero as $n \to \infty$.
By the assumptions regarding the functions $v$ and $q$ and an application of  the mean value theorem this is equivalent to
the condition
\begin{equation}\label{hd8}
	\sup_{f \in \Theta(\beta,L)} \sum_{i=1}^{n} \Big ( f(t_{i,n})  - n \int_{t_{i-1,n}}^{t_{i,n}} f(s)\dd s \Big )^2 \to 0.
\end{equation}
Let $f\in \Theta(\beta,L)$ with Fourier expansion $f(\cdotb) = \sum_{k \in \Z} \theta_k \eb_k(\cdotb)$.
For any $K \in \N$ (the appropriate value of $K=K(n)$ for our purposes will be specified below) we define 
the functions 
\begin{align*}
f_K(\cdotb) & = \sum_{\abs{k} \leq K} \theta_k \eb_k(\cdotb)
\text{~~~~
and } \\
f_K^\top(\cdotb) & = f(\cdotb) - f_K(\cdotb) = \sum_{\abs{k} > K} \theta_k \eb_k(\cdotb),
\end{align*}
respectively.
Consequently,
\begin{align}    \label{hd7}
	\sum_{i=1}^n \Big (f(t_{i,n})  - n \int_{t_{i-1,n}}^{t_{i,n}} f(s)\dd s \Big )^2 \leq 3 \sum_{i=1}^n A_{i,n}^2 + 3 \sum_{i=1}^n B_{i,n}^2 + 3 \sum_{i=1}^n C_{i,n}^2,
\end{align}
where
\begin{align*}
	A_{i,n} &= f_K(t_{i,n})  - n \int_{t_{i-1,n}}^{t_{i,n}} f_K(s)\dd s ,\\
	B_{i,n} &= f_K^\top(t_{i,n}), \quad \text{and}\\
	C_{i,n} &= n \int_{t_{i-1,n}}^{t_{i,n}} f_K^\top(s)\dd s.
\end{align*}
We now show for $K=n$ that the estimates
\begin{align}
  \sum_{i=1}^n A_{i,n}^2 &=  \Oc (\max \{ n^{-1},n^{1-2\beta} \} ) , \label{hd4}
    \\
  \sum_{i=1}^n B_{i,n}^2 &=  \Oc ( n^{1-2\beta}  )   ,     \label{hd5}
    \\
   \sum_{i=1}^n C_{i,n}^2 &= \Oc ( n^{1-2\beta} )   ,     \label{hd6}
\end{align}
hold uniformly with respect to $f \in \Theta(\beta,L)$. Then, assertion \eqref{hd8} follows from \eqref{hd7} and the assumption $\beta >1/2$.
\medskip

\noindent
{\it Proof of \eqref{hd4}:} 
 The quantity
\begin{equation*}
	\frac{1}{n} \sum_{j=1}^n \lvert A_{j,n} \rvert^2  
\end{equation*}
can be interpreted as the (average) energy of the discrete signal $A^{(n)} = (A_{1,n},\ldots,A_{n,n})$.
Define 
\begin{equation*}
	F_j = \frac{1}{n} \sum_{k=1}^n A_{k,n} e^{-2\pi \ii k j/n}
\end{equation*}
as  the discrete Fourier transform of the signal $A^{(n)}$, then 
Parseval's identity for the discrete Fourier transform yields
\begin{equation*}
	\frac{1}{n} \sum_{j=1}^n \lvert A_{j,n} \rvert^2 = \sum_{j=1}^n \lvert F_j \rvert^2,
\end{equation*}
and we have to derive an estimate for $n \sum_{j=1}^n \lvert F_j \rvert^2   $.
For this purpose, we recall the notation of  $A_{j,n}$  and note that
\begin{align*}
	F_j &= \frac{1}{n} \sum_{k=1}^{n} \Big  (f_K(t_{k,n}) - n \int_{t_{k-1}}^{t_k} f_K(s) \dd s \Big ) e^{-2\pi \ii kj/n}\\
	&= \frac{1}{n} \sum_{k=1}^{n} \left( \sum_{1 \leq \lvert l \rvert \leq K} \theta_l e^{-2\pi \ii l k/n} - n \int_{t_{k-1}}^{t_k} \sum_{1 \leq \lvert l \rvert \leq K} \theta_l e^{-2\pi \ii ls} \dd s \right)  e^{-2\pi \ii kj/n}\\
	&= \frac{1}{n} \sum_{k=1}^{n} \left( \sum_{1 \leq \lvert l \rvert \leq K} \theta_l e^{-2\pi \ii l k/n} + n \sum_{1 \leq \lvert l \rvert \leq K} \frac{\theta_l}{2\pi \ii l} \left[ e^{-2\pi \ii lk/n} - e^{-2\pi \ii l(k-1)/n} \right]   \right)  e^{-2\pi \ii kj/n}\\
	&= \frac{1}{n} \sum_{k=1}^{n} \left( \sum_{1 \leq \lvert l \rvert \leq K} \theta_l e^{-2\pi \ii l k/n} + n \sum_{1 \leq \lvert l \rvert \leq K} \frac{\theta_l}{2\pi \ii l} \left[ 1 - e^{2\pi \ii l/n} \right] e^{-2\pi \ii lk/n} \right)  e^{-2\pi \ii kj/n}\\
	&= \frac{1}{n} \sum_{k=1}^{n} \sum_{1 \leq \lvert l \rvert \leq K} \theta_l \left[ 1 + \frac{n}{2\pi \ii l} ( 1- e^{2\pi \ii l/n} )  \right] e^{-2\pi \ii kl/n} e^{-2\pi \ii kj /n}. 
\end{align*}
From now on, we take $K=n$ and write
\begin{align*}
	F_j^+ &= \frac{1}{n} \sum_{k=1}^{n} \sum_{l=1}^n \theta_l \left[ 1 + \frac{n}{2\pi \ii l} ( 1- e^{2\pi \ii l/n} )  \right] e^{-2\pi \ii kl/n} e^{-2\pi \ii kj /n}, \quad \text{and}\\
	F_j^- &= \frac{1}{n} \sum_{k=1}^{n} \sum_{l=-n}^{-1} \theta_l \left[ 1 + \frac{n}{2\pi \ii l} ( 1- e^{2\pi \ii l/n} )  \right] e^{-2\pi \ii kl/n} e^{-2\pi \ii kj /n}.
\end{align*}
Since $\sum_{j=1}^n \lvert F_j \rvert^2 \leq 2 \sum_{j=1}^n \lvert F_j^+ \rvert^2 + 2\sum_{j=1}^n \lvert F_j^- \rvert^2$, it is sufficient to consider $\sum_{j=1}^n \lvert F_j^+ \rvert^2$ (the term involving $F_j^-$ is treated analogously).
We have
\begin{align*}
	F_j^+ &= \frac{1}{n} \sum_{k=1}^{n} \sum_{l=1}^n \theta_l \left[ 1 + \frac{n}{2\pi \ii l} ( 1- e^{2\pi \ii l/n} )  \right] e^{-2\pi \ii kl/n} e^{-2\pi \ii kj /n}\\
	&= \frac{1}{n} \sum_{l=1}^n \theta_l  \left[ 1 + \frac{n}{2\pi \ii l} ( 1- e^{2\pi \ii l/n} )  \right]  \sum_{k=1}^{n}  e^{-2\pi \ii kl/n} e^{-2\pi \ii kj /n}\\
	&= \theta_{l(j)} \left[ 1 + \frac{n}{2\pi \ii l(j)} ( 1- e^{2\pi \ii l(j)/n} )  \right] 
\end{align*}
where $l(j) = n - j$ for $j=1,\ldots,n-1$ and $l(n) = n$.
Here, we used the well-known fact that for any integer $m \in \Z$
\begin{equation}\label{eq:geom_ser_root}
	\sum_{k=1}^n e^{-2\pi \ii km/n} = \begin{cases}
		n, & \text{if } m \in n\Z,\\
		0, & \text{if } m \notin n \Z.
	\end{cases}  
\end{equation}
Thus, we obtain (uniformly with respect to $\Theta$) 
\begin{align*}
	\sum_{j=1}^n \lvert F_j^+ \rvert^2 &=  \sum_{j=1}^n \lvert \theta_j \rvert^2 \left|  1 + \frac{n}{2\pi \ii j} ( 1- e^{2\pi \ii j/n} )  \right|^2\\
	&=  \sum_{j=1}^n \lvert \theta_j \rvert^2 \left| \frac{n}{2\pi \ii j} ( e^{2\pi \ii j/n} - 1 - 2\pi \ii j/n ) \right|^2\\
	&\asymp n^2\sum_{j=1}^n \frac{\lvert \theta_j \rvert^2}{\lvert j \rvert^2} \left|  e^{2\pi \ii j/n} - 1 - 2\pi \ii j/n \right|^2\\
	&\lesssim n^2\sum_{j=1}^n \frac{\lvert \theta_j \rvert^2}{\lvert j \rvert^2} \cdot \left|  j/n \right|^4
	= n^{-2} \sum_{j=1}^n \lvert \theta_j \rvert^2 \lvert j \rvert^2\\
	&= n^{-2} \sum_{j=1}^n \lvert \theta_j \rvert^2 \lvert j \rvert^{2\beta} \lvert j \rvert^{2-2\beta}
	\leq n^{-2} L^2 \max \{1,n^{2-2\beta}\} \\ 
	&\lesssim \max \{ n^{-2},n^{-2\beta} \}.
\end{align*}
An analogous argument for the term 
 $\sum_{j=1}^n \lvert F_j^- \rvert^2  $ proves  the estimate \eqref{hd4}.

\medskip

\noindent
{\it Proof of \eqref{hd5}:} We have
\begin{align*}
	\sum_{j=1}^n \lvert f_K^\top(t_{j,n}) \rvert^2 &= \sum_{j=1}^n \left|  \sum_{\abs{k} > K} \theta_k \exp(-2\pi \ii kj/n) \right|^2\\
	&\leq 2\sum_{j=1}^n \left|  \sum_{k > K} \theta_k \exp(-2\pi \ii kj/n) \right|^2 + 2\sum_{j=1}^n \left|  \sum_{k < -K} \theta_k \exp(-2\pi \ii kj/n) \right|^2,
\end{align*}
and it is again sufficient to consider the sum running over $k>K$.
Using \eqref{eq:geom_ser_root} again, we get
\begin{align*}
	\sum_{j=1}^n \left|  \sum_{k > K} \theta_k \exp(-2\pi \ii kj/n) \right|^2&= \sum_{k,l > K} \theta_k \overline{\theta_{l}}  \sum_{j=1}^n \exp(-2\pi \ii (k-l)j/n)\\
	&= n \sum_{\substack{k,l > K\\ k-l \in n\Z}} \theta_k \overline{\theta_{l}}.
\end{align*}
Taking $K = n$ here as well yields
\begin{align*}
	\sum_{j=1}^n \left|  \sum_{k > K} \theta_k \exp(-2\pi \ii kj/n) \right|^2 &= \sum_{m=n+1}^{2n} \left| \sum_{r=0}^\infty \theta_{m + rn} \right|^2.
\end{align*}
Now,
\begin{align*}
	\left| \sum_{r=0}^\infty \theta_{m + rn} \right|^2 &\leq \left(  \sum_{r=0}^\infty \lvert \theta_{m + rn}\rvert^2 (m+rn)^{2\beta} \right) \left( \sum_{r=0}^\infty (m+rn)^{-2\beta} \right),
\end{align*}
and  we obtain 
\begin{align*}
	\sum_{m=n+1}^{2n} \left| \sum_{r=0}^\infty \theta_{m + rn} \right|^2 &\leq n L^2 \sum_{r=0}^\infty (n+rn)^{-2\beta}
	\lesssim n^{1-2\beta}
\end{align*}
 uniformly over $f \in \Theta(\beta,L)$, which establishes \eqref{hd5}.
\medskip

\noindent
{\it Proof of \eqref{hd6}:}
Using Jensen's inequality and Parseval's identity we obtain
\begin{align*}
	\sum_{i=1}^n \left( n \int_{t_{i-1,n}}^{t_{i,n}} f_K^\top(s)\dd s\right) ^2 &\leq n \int_0^1 (f_K^\top(s))^2 \dd s= n \sum_{\abs{k} > K} \lvert \theta_k \rvert^2\leq nL^2 K^{-2\beta}
\end{align*}
 uniformly over $f \in \Theta(\beta,L)$,
which is of order $n^{1-2\beta}$ if we choose $K=n$.

\medskip

\noindent\underline{Verification of condition  \ref{triC2}}:
We have to show that
\begin{equation}\label{eq:def:Dn}
D_n:= 	\min_{\alphab_n \in \R^n} \lVert F_f(\cdotb) - \sum_{j=1}^n \alpha_j K_\Xi(\cdotb,t_{j,n}) \rVert_{\RKHS(\Xi)}^2 = o(n^{-1})
\end{equation}
uniformly over all $f \in \Theta(\beta,L)$.
Via the isomorphism $\psi$ introduced in Equation~\eqref{psi} we have for $g=\psi F_f$
\begin{align*}
D_n	&= \min_{\alphab_n \in \R^n} \lVert g(\cdotb) - \sum_{j=1}^n \alpha_j v(t_{j,n}) \1_{[0,q(t_{j,n})]}(\cdotb)  \rVert_{L^2([0,T])}^2\\
	&= \min_{\alphab_n \in \R^n} \lVert g(\cdotb) - \sum_{j=1}^n \alpha_j \1_{(q(t_{j-1,n}),q(t_{j,n})]}(\cdotb) \rVert^2_{L^2([0,T])}\\
	&=\min_{\alphab_n \in \R^n} \int_0^T (g(u) - \sum_{j=1}^n \alpha_j \1_{(q(t_{j-1,n}),q(t_{j,n})]}(u) )^2 \dd u\\
	&= \min_{\alphab_n \in \R^n} \int_0^1 (g(q(w)) - \sum_{j=1}^n \alpha_j \1_{(q(t_{j-1,n}),q(t_{j,n})]}(q(w)))^2 q'(w) \dd w\\
	&= \min_{\alphab_n \in \R^n} \int_0^1 (g(q(w)) - \sum_{j=1}^n \alpha_j \1_{(t_{j-1,n},t_{j,n}]}(w))^2 q'(w) \dd w.
\end{align*}
Assuming that $q'$ is bounded from above we obtain
\begin{align*}
D_n	\leq C(q) &\min_{\alphab_n \in \R^n} \int_0^1 (g(q(w)) - \sum_{j=1}^n \alpha_j \1_{(t_{j-1,n},t_{j,n}]}(w))^2 \dd w.
\end{align*}
Note that 
\begin{equation*}
	g(q(w)) = g_1(q(w)) - g_2(q(w))
\end{equation*}
with
\begin{align*}
	g_1(q(w)) = \frac{f(w)}{v(w)q'(w)},\quad \text{and}\quad g_2(q(w)) = \frac{F_f(w)v'(w)}{v^2(w)q'(w)}.
\end{align*}
With $f_n(\cdotb) = \sum_{\abs{k} \leq n} \theta_k \eb_k(\cdotb)$ (and $f_n^\top(\cdotb) = f(\cdotb) - f_n(\cdotb)$) we define
\begin{align*}
	\alpha_j^{(1)} = \frac{f_n(t_{j,n})}{v(t_{j,n})q'(t_{j,n})}
	\quad \text{ and } \quad
	\alpha_j^{(2)} = \frac{F_{f}(t_{j,n})v'(t_{j,n})}{v^2(t_{j,n})q'(t_{j,n})}.
\end{align*}
Using these notations we get
\begin{equation}\label{eq:bound:Dn}
D_n = 	\min_{\alphab_n \in \R^n} \lVert F_f(\cdotb) - \sum_{j=1}^n \alpha_j K_\Xi(\cdotb,t_{j,n}) \rVert_{\RKHS(\Xi)}^2 \lesssim I_1 + I_2
\end{equation}
where
\begin{align} \label{hd11} 
	I_1 = \int_0^1 \Big ( g_1(q(w)) - \sum_{j=1}^{n} \alpha_j^{(1)} \1_{[t_{j-1,n},t_{j,n})}(w)\Big )^2 \dd w , \\
\label{hd12} 
	I_2 = \int_0^1 \Big (g_2(q(w)) - \sum_{j=1}^{n} \alpha_j^{(2)} \1_{[t_{j-1,n},t_{j,n})}(w) \Big )^2 \dd w.
\end{align}
We investigate the two terms $I_1$ and $I_2$ separately.

\medskip

\noindent\emph{Bound for $I_1$}: We use the estimate 
\begin{align}\label{eq:I_1}
	I_1 &\lesssim I_{11} + I_{12}, \end{align}
where
\begin{align*}
	I_{11}  &= \int_0^1 \left( \frac{f_n(w)}{v(w)q'(w)} - \sum_{j=1}^{n} \alpha_j^{(1)} \1_{[t_{j-1,n},t_{j,n})}(w)\right) ^2 \dd w , \\
		I_{12}  &=  \int_0^1 \left( \frac{f_n^\top(w)}{v(w)q'(w)}\right) ^2 \dd w.
\end{align*}
For the first integral $I_{11}$ on the right-hand side of \eqref{eq:I_1}, we have
\begin{align}
I_{11} 	&= \sum_{j=1}^n \int_{t_{j-1,n}}^{t_{j,n}} \left( \frac{f_n(w)}{v(w)q'(w)} - \frac{f_n(t_{j,n})}{v(t_{j,n})q'(t_{j,n})} \right) ^2 \dd w\notag\\
	&\lesssim \sum_{j=1}^n \int_{t_{j-1,n}}^{t_{j,n}} \left( \frac{f_n(w)}{v(w)q'(w)} - \frac{f_n(t_{j,n})}{v(w)q'(w)} \right) ^2 \dd w\label{eq:I1_1}\\
	&+ \sum_{j=1}^n \int_{t_{j-1,n}}^{t_{j,n}} \left( \frac{f_n(t_{j,n})}{v(w)q'(w)} - \frac{f_n(t_{j,n})}{v(t_{j,n})q'(t_{j,n})} \right) ^2 \dd w.\label{eq:I1_2}
\end{align}
First, we further decompose \eqref{eq:I1_1} as
\begin{align*}
	\sum_{j=1}^n \int_{t_{j-1,n}}^{t_{j,n}} \left( \frac{f_n(w)}{v(w)q'(w)} - \frac{f_n(t_{j,n})}{v(w)q'(w)} \right) ^2 \dd w &\lesssim C(v,q) \sum_{j=1}^n \int_{t_{j-1,n}}^{t_{j,n}} ( f_n(w) - f_n(t_{j,n})) ^2 \dd w\\
	&\hspace{-13em}\lesssim \sum_{j=1}^n \int_{t_{j-1,n}}^{t_{j,n}} \lvert f_n^+(w) - f_n^+(t_{j,n}) \rvert ^2 \dd w + \sum_{j=1}^n \int_{t_{j-1,n}}^{t_{j,n}} \lvert f_n^-(w) - f_n^-(t_{j,n})\rvert ^2 \dd w~, 
\end{align*}
where
\begin{align*}
	f_n^+ = \sum_{k=1}^n \theta_k \eb_k(\cdotb) \quad \text{and} \quad f_n^- =  \sum_{k=-n}^{-1} \theta_k \eb_k(\cdotb).
\end{align*}
In the sequel, we consider only the term involving $f_n^+$ since the sum involving $f_n^-$ can be bounded using the same argument.
We have the identity
\begin{align*}
	\lvert f_n^+(w) - f_n^+(t_{j,n}) \rvert^2 &= \lvert \sum_{k=1}^n \theta_k (\eb_k(w) - \eb_k(t_{j,n})) \rvert^2\\
	&\hspace{-5em}= \sum_{k,l = 1}^n \theta_k \overline{\theta}_{l} \exp(-2\pi \ii kj/n)[\exp(2\pi \ii k(j/n - w)) - 1]  \\
& \hspace{5em}	\cdot \exp(2\pi \ii l j/n) [\exp(-2\pi \ii l(j/n - w)) - 1].
\end{align*}
From this identity we obtain (exploiting \eqref{eq:geom_ser_root} again)
\begin{align*}
	&\sum_{j=1}^n \int_{t_{j-1,n}}^{t_{j,n}} \lvert f_n^+(w) - f_n^+(t_{j,n})\rvert ^2 \dd w\\
	&\hspace{1em}= \sum_{j=1}^n \sum_{k,l=1} \theta_k \overline{\theta}_l \exp(2\pi \ii (l-k)j/n) \int_{\frac{j-1}{n}}^{\frac{j}{n}} [\exp(2\pi \ii k(j/n - w)) - 1] \\
	& \hspace{7.5cm} \cdot [\exp(-2\pi \ii l(j/n - w)) - 1] \dd w\\
	&\hspace{1em}= \sum_{j=1}^n \sum_{k,l=1} \theta_k \overline{\theta}_l \exp(2\pi \ii (l-k)j/n) \int_{0}^{\frac{1}{n}} [\exp(2\pi \ii k(1/n - w)) - 1]\\
 	& \hspace{7.55cm} \cdot [\exp(-2\pi \ii l(1/n - w)) - 1] \dd w\\
	&\hspace{1em}= n \sum_{k=1}^n \abs{\theta_k}^2 \int_{0}^{\frac{1}{n}} [\exp(2\pi \ii k(1/n - w)) - 1]   [\exp(-2\pi \ii k(1/n - w)) - 1] \dd w\\
	&\hspace{1em}\leq C n \sum_{k=1}^n \abs{\theta_k}^2 k^2 n^{-3}\\
	&\hspace{1em}\leq C n^{-2} \sum_{k=1}^n \abs{\theta_k}^2 k^{2\beta} k^{-2\beta + 2}\\
	&\hspace{1em}\leq C(L)n^{-2} \max \{ 1, n^{-2\beta+2} \} \lesssim \max \{ n^{-2}, n^{-2\beta} \} = o(n^{-1}).
\end{align*}
\noindent 
To derive an estimate of \eqref{eq:I1_2} we note that for 
  any $n \in \N$,
\begin{align*}
	\lVert f_n \rVert^2_\infty =  \sup_{x \in [0,1]} \lvert \sum_{\abs{k} \leq n} \theta_k \eb_k(x) \rvert^2 &\leq \left( \sum_{\abs{k} \leq n} \lvert \theta_k \rvert^2 k^{2\beta} \right) \left( \sum_{\abs{k} \leq n} k^{-2\beta} \right) \leq C(L,\beta) 
\end{align*}
(the same estimate holding true for $f$ instead of $f_n$ which formally corresponds to $n = \infty$).
Hence,
\begin{align*}
	\sum_{j=1}^n &\int_{t_{j-1,n}}^{t_{j,n}} \left( \frac{f_n(t_{j,n})}{v(w)q'(w)} - \frac{f_n(t_{j,n})}{v(t_{j,n})q'(t_{j,n})} \right) ^2 \dd w\\
	&= \sum_{j=1}^n \lvert f_n(t_{j,n}) \rvert^2 \int_{t_{j-1,n}}^{t_{j,n}} \left( \frac{1}{v(w)q'(w)} - \frac{1}{v(t_{j,n})q'(t_{j,n})} \right) ^2 \dd w\\
	&\leq C(L,\beta) \sum_{j=1}^n \int_{t_{j-1,n}}^{t_{j,n}} \left|  \frac{v(t_{j,n})q'(t_{j,n}) - v(w)q'(w)}{v(w)q'(w)v(t_{j,n})q'(t_{j,n})} \right|^2 \dd w\\
	&\leq C(L,\beta,v,q) \sum_{j=1}^n \int_{t_{j-1,n}}^{t_{j,n}} \lvert v(t_{j,n})q'(t_{j,n}) - v(w)q'(w) \rvert^2 \dd w.
\end{align*}
Because the product of a $\gamma_1$-Hölder function and a $\gamma_2$-Hölder function is (at least) Hölder with index $\min\{ \gamma_1,\gamma_2 \}$ we obtain from our assumptions that the function $vq'$ is Hölder continuous with index $\gamma>1/2$.
Thus,
\begin{align*}
	\sum_{j=1}^n \int_{t_{j-1,n}}^{t_{j,n}} \left( \frac{f_n(t_{j,n})}{v(w)q'(w)} - \frac{f_n(t_{j,n})}{v(t_{j,n})q'(t_{j,n})} \right) ^2 \dd w &\leq C(L,\beta,v,q) \sum_{j=1}^n \int_{t_{j-1,n}}^{t_{j,n}} n^{-2\gamma} \dd w\\
	&= C(L,\beta,v,q) n^{-2\gamma},
\end{align*}
and this is $o(n^{-1})$ since $\gamma > 1/2$. Combining these arguments  we obtain 
\begin{equation*}
    I_{11} = o(n^{-1}). 
\end{equation*}

\medskip

\noindent Finally, the second integral $I_{12}$ on the right-hand side of \eqref{eq:I_1} can bounded as follows: 
\begin{align*}
I_{12} &= 	\int_0^1 \left( \frac{f_n^\top(w)}{v(w)q'(w)}\right) ^2 \dd w \leq C(v,q) \int_0^1 \lvert f_n^\top(w)\rvert^2 \dd w\\
	&\leq C(v,q) \sum_{\lvert k \rvert > n} \lvert \theta_k \rvert^2
	\leq C(v,q)L^2 n^{-2\beta} = o(n^{-1}).
\end{align*}
Observing the estimate \eqref{hd8} we finally obtain  $I_1= o(n^{-1})$.
\smallskip 

\noindent\emph{Bound for $I_2$}: In analogy to the decomposition of the term $I_{11}$ on the right-hand side of \eqref{eq:I_1} we have
\begin{align}
	I_2 &= \int_0^1 \left( \frac{F_{f}(w)v'(w)}{v^2(w)q'(w)} - \sum_{j=1}^{n} \alpha_j^{(2)} \1_{[t_{j-1,n},t_{j,n})}(w)\right) ^2 \dd w\notag\\
	&\lesssim \sum_{j=1}^n \int_{t_{j-1,n}}^{t_{j,n}} \left( \frac{F_{f}(w)v'(w)}{v^2(w)q'(w)} - \frac{F_{f}(t_{j,n})v'(t_{j,n})}{v^2(w)q'(w)} \right) ^2 \dd w\label{eq:I_2:1}\\
	&\hspace*{1em} + \sum_{j=1}^n \int_{t_{j-1,n}}^{t_{j,n}} \left( \frac{F_{f}(t_{j,n})v'(t_{j,n})}{v^2(w)q'(w)} - \frac{F_{f}(t_{j,n})v'(t_{j,n})}{v^2(t_{j,n})q'(t_{j,n})} \right) ^2 \dd w\label{eq:I_2:2}.
\end{align}

Note that $F_f$ is Lipschitz since it is continuously differentiable (recall that $f$ itself is continuous since $\beta > 1/2$) and $v'$ is Hölder with index $\gamma>1/2$ due to our assumptions.
Thus, the term \eqref{eq:I_2:1} can be bounded as
\begin{align*}
	\sum_{j=1}^n \int_{t_{j-1,n}}^{t_{j,n}} \left( \frac{F_{f}(w)v'(w)}{v^2(w)q'(w)} - \frac{F_{f}(t_{j,n})v'(t_{j,n})}{v^2(w)q'(w)} \right) ^2 \dd w\leq C(v,q) n^{-2\gamma}
\end{align*}
which is of order $o(n^{-1})$.
For the term \eqref{eq:I_2:2} we obtain using our assumptions that
\begin{align*}
	&\sum_{j=1}^n \int_{t_{j-1,n}}^{t_{j,n}} \left( \frac{F_{f}(t_{j,n})v'(t_{j,n})}{v^2(w)q'(w)} - \frac{F_{f}(t_{j,n})v'(t_{j,n})}{v^2(t_{j,n})q'(t_{j,n})} \right) ^2 \dd w\\ &\hspace{-2em}\leq C(\beta,L,v,q) \sum_{j=1}^n \int_{t_{j-1,n}}^{t_{j,n}} (v^2(t_{j,n})q'(t_{j,n}) - v^2(w)q'(w))^2 \dd w.
\end{align*}
Using the same arguments as for the bound of \eqref{eq:I1_2}, this term can be shown to be of order $o(n^{-1})$.
Since both terms $I_1$ and $I_2$ are of order $o(n^{-1})$ the assertion \eqref{eq:def:Dn} follows from \eqref{eq:bound:Dn}.

\subsection{Proof of Theorem~\ref{thm:hoelder}}\label{app:hoelder}

As in the proof of Theorem~\ref{thm:sobolev}, we have to verify the the two conditions \ref{triC1} and \ref{triC2} from Theorem~\ref{thm:general:triangular}.
\medskip

\noindent\underline{Verification of condition \ref{triC1}}:
As in the Sobolev case it is sufficient to show that
\begin{equation*}
	\sup_{f \in \Fc(\alpha, L, M) } \sum_{i=1}^{n} \Big  ( f(t_{i,n})  - n \int_{t_{i-1,n}}^{t_{i,n}} f(s)\dd s \Big )^2 \to 0.
\end{equation*}
By the mean value theorem $n \int_{t_{i-1,n}}^{t_{i,n}} f(s)\dd s = f(\zeta_{i,n})$ for some $t_{i-1,n} \leq \zeta_{i,n} \leq t_{i,n}$.
Thus, since $f \in \Fc(\alpha,L,M)$,
\begin{align*}
    \sum_{i=1}^{n} \Big ( f(t_{i,n})  - n \int_{t_{i-1,n}}^{t_{i,n}} f(s)\dd s \Big )^2 &= \sum_{i=1}^n (f(t_{i,n}) - f(\zeta_{i,n}))^2\\
    &\leq \sum_{i=1}^n L^2 \lvert t_{i,n} - \zeta_{i,n} \rvert^2\\
    &\leq L^2 n^{-2\alpha +1},
\end{align*}
and the last term converges to zero uniformly over $f \in \Fc(\alpha,L,M)$ whenever $\alpha > 1/2$.
\medskip 

\noindent\underline{Verification of  condition \ref{triC2}}: The proof is based on nearly the same reduction as in the Sobolev case.
Again, we consider the bound
\begin{equation}\label{eq:H:I1I2}
	\min_{\alphab_n \in \R^n} \lVert F_f(\cdotb) - \sum_{j=1}^n \alpha_j K_\Xi(\cdotb,t_{j,n}) \rVert_{\RKHS(\Xi)}^2 \lesssim I_1 + I_2
\end{equation}
where we define $I_1$ and $I_2$ as in Appendix~\ref{app:sobolev} (see the equations \eqref{hd11} and \eqref{hd12}) with the only exception that we now put
\begin{equation*}
    \alpha_j^{(1)} = \frac{f(t_{j,n})}{v(t_{j,n})q'(t_{j,n})}
\end{equation*}
in the definition of $I_1$.
Then,
\begin{align*}
	I_1 &\lesssim \int_0^1 \left( \frac{f(w)}{v(w)q'(w)} - \sum_{j=1}^{n} \alpha_j^{(1)} \1_{[t_{j-1,n},t_{j,n})}(w)\right) ^2 \dd w\\
	&= \sum_{j=1}^n \int_{t_{j-1,n}}^{t_{j,n}} \left( \frac{f(w)}{v(w)q'(w)} - \frac{f(t_{j,n})}{v(t_{j,n})q'(t_{j,n})} \right) ^2 \dd w\\
	&\lesssim \sum_{j=1}^n \int_{t_{j-1,n}}^{t_{j,n}} \left( \frac{f(w)}{v(w)q'(w)} - \frac{f(t_{j,n})}{v(w)q'(w)} \right) ^2 \dd w\\
	&\hspace{1em}+ \sum_{j=1}^n \int_{t_{j-1,n}}^{t_{j,n}} \left( \frac{f(t_{j,n})}{v(w)q'(w)} - \frac{f(t_{j,n})}{v(t_{j,n})q'(t_{j,n})} \right) ^2 \dd w.
\end{align*}
The first integral can be bounded as
\begin{align*}
    \sum_{j=1}^n \int_{t_{j-1,n}}^{t_{j,n}} \left( \frac{f(w)}{v(w)q'(w)} - \frac{f(t_{j,n})}{v(w)q'(w)} \right) ^2 \dd w &\leq C(v,q) \sum_{j=1}^n \int_{t_{j-1,n}}^{t_{j,n}} (f(w) - f(t_{j,n}))^2 \dd w\\
    &\leq C(v,q) L^2 n^{-2\alpha},
\end{align*}
which converges to zero as $n$ increases.
The second integral can be bounded as in the Sobolev case using the assumption that $f$ is bounded (as one can easily see, the assumption of uniform boundedness can be dropped if $v \cdot q'$ is constant; this is for instance satisfied in the case of Brownian motion).
Hence, $I_1$ converges to $0$.

The term $I_2$ in \eqref{eq:H:I1I2} can be bounded exactly as the corresponding term in the Sobolev ellipsoid case.
This finishes the proof of the theorem.

\subsection{Proof of Theorem~\ref{thm:bridge}}

To prepare the proof, we recall an alternative (but equivalent) characterization of asymptotic equivalence in the framework of statistical decision theory.
Let $\Expm = (\Xc,\Xs,(\Pb_{\theta})_{\theta \in \Theta})$ be a statistical experiment.
In decision theory one considers a decision space $(\Ac,\As)$ where the set $\Ac$ contains the potential decisions (or actions) that are at the observers disposal and $\As$ is a $\sigma$-field on $\Ac$.
In addition, there is a \emph{loss function}
\begin{equation*}
\ell \colon \Theta \times \Ac \to [0,\infty), \qquad (\theta,a) \mapsto \ell(\theta,a)
\end{equation*}
with the interpretation that a loss $\ell(\theta,a)$ occurs if the statistician chooses the action $a \in {\cal A} $ and $\theta \in \Theta $ is the true state of nature.
A (randomized) decision rule is a Markov kernel $\rho \colon \Xc \times \As \to [0,1]$, and the associated \emph{risk} is
\begin{equation*}
R_\theta(\Expm,\rho,\ell) = \int_\Xc \left( \int_\Ac \ell(\theta, a) \rho(x,\dd a) \right) \Pb_\theta(\dd x).
\end{equation*}
Then, the deficiency between two experiments $\Expm_1$ and $ \Expm_2 $ is exactly the quantity
\begin{equation}\label{eq:deficiency:loss}
\delta(\Expm_1,\Expm_2) = \inf_{\rho_1} \sup_{\rho_2} \sup_\theta \sup_\ell \lvert R_\theta(\Expm_1,\rho_1,\ell) - R_\theta(\Expm_2,\rho_2,\ell) \rvert
\end{equation}
(see \cite{mariucci2016lecam}, Theorem~2.7, and the references cited there), where the supremum is taken over all loss functions $\ell$ with $0 \leq \ell(\theta,a) \leq 1$ for all $\theta \in \Theta$ and $a \in \Ac$, all admissible parameters $\theta \in \Theta$, and decision rules $\rho_2$ in the second experiment.
The infimum is taken over all decision rules $\rho_1$ in the first experiment.

After these preliminaries, let us go on to the proof of the theorem.
We consider the decision space $(\Ac,\As) = (\R,\Bs(\R))$ and the loss function 
\begin{equation*}
\ell \colon \Theta(\beta,L) \times \Ac \to \{0,1\}, \quad (f,a) \mapsto \ell(f,a) = \begin{cases}
1, & \text{if } \int_0^1 f(x)\dd x \neq a,\\
0, & \text{if } \int_0^1 f(x)\dd x = a.
\end{cases}
\end{equation*}
In the experiment $\Expm_{2,n}$ we observe the whole path $Y = \{ Y_t, \, t \in [0,1] \}$ satisfying
\begin{equation*}
Y_t = \int_0^t f(s) \dd s + \frac{1}{\sqrt{n}} B_t, \quad t \in [0,1]~,
\end{equation*}
where $(B_t)_{t \in [0,1]}$ is a Brownian Bridge, and we consider the (non-randomized) decision rule $\rho_2$ defined by
\begin{equation*}
\rho_2(h) = h(1) - h(0)~,~~~~ h \in {\cal C}([0,1], \R) .
\end{equation*}
This directly yields $\rho_2(Y) = \int_0^1 f(s)\dd s$ since $\Xi_0 = \Xi_1 = 0$ for the Brownian bridge.
Hence,
\begin{equation*}
R_f(\Expm_{2,n},\rho_2,\ell) = \int \ell(f,\rho_2(Y))  \Pb^f_{2,n}(\dd Y) = 0
\end{equation*}
for all $f \in \Theta(\beta,L)$ and $n \in \N$.
From \eqref{eq:deficiency:loss} we thus obtain
\begin{align*}
\delta(\Expm_{1,n},\Expm_{2,n}) &\geq \inf_{\rho_1} \sup_{f \in \Theta(\beta,L)} \lvert R_f(\Expm_{1,n},\rho_1,\ell) \rvert.
\end{align*}
Recall the notation $\eb_n(x) = \exp(-2\pi \ii n x)$
and introduce the functions
$f_0 \equiv 0$ and
\begin{equation*}
f_n(x) = \sqrt{\frac{2}{3}}\frac{L}{n^\beta} \left[ 1 - \frac{1}{2} \eb_n(x) - \frac{1}{2} \eb_{-n}(x)  \right]
= \sqrt{\frac{2}{3}}\frac{L}{n^\beta} \left[ 1 - \frac{1}{2} \cos( 2\pi n x)  \right]
\end{equation*}
for $n \in \N$.
It is easily seen that $f_n$ belongs to $\Theta(\beta,L)$. 
Note that by construction $f_n(j/n) = 0$ for $j=1,\ldots,n$ and thus in the experiment $\Expm_{1,n}$ the identity $\Pb_{1,n}^{f_0} = \Pb_{1,n}^{f_n}$ holds.
For the considered loss function we have
\begin{align*}
R_f(\Expm_{1,n},\rho_1,\ell) &= \int_{\R^n} 
 \rho_1(\Yb_n) \left( \R \backslash \left\{ \int_0^1 f(x)\dd x \right\} \right) \Pb_{1,n}^f(\dd \Yb_n),
\end{align*}
where $\rho_1$ is any (potentially randomized) decision rule.
Because 
$$\int_0^1 f_0(x)\dd x= 0 \neq \sqrt{2/3}Ln^{-\beta} = \int_0^1 f_n(x)\dd x
$$
at least one of the quantities $\rho_1(\Yb_n) ( \R \backslash \{ \int_0^1 f_0(x)\dd x \} )$ and $\rho_1(\Yb_n) ( \R \backslash \{ \int_0^1 f_n(x)\dd x \} )$  must be $\geq 1/2$ for any $\Yb_n$ (otherwise there is a contradiction).
Thus, setting $A_\diamond = \{ \Yb_n \in \R^n : \rho_1(\Yb_n) ( \R \backslash \{ \int_0^1 f_\diamond(x)\dd x \} ) \geq 1/2 \}$ for $\diamond \in \{ 0, n \}$ one has $A_0 \cup A_n = \R^n$.
As a consequence, either $\Pb_{1,n}^{f_0}(A_0) = \Pb_{1,n}^{f_n}(A_0) \geq 1/2$ or $\Pb_{1,n}^{f_0}(A_n) = \Pb_{1,n}^{f_n}(A_n) \geq 1/2$ holds.
Without loss of generality, we assume that $\Pb_{1,n}^{f_0}(A_0) = \Pb_{1,n}^{f_n}(A_0) \geq 1/2$ (the other case follows by exactly the same argument).
In this case, using the definition of the set $A_0$,
\begin{align*}
\delta(\Expm_{1,n},\Expm_{2,n}) &\geq \inf_{\rho_1} \sup_{f \in \{ f_0,f_n \}} R_f(\Expm_{1,n},\rho_1,\ell) \\
&\geq \inf_{\rho_1} \int_{\R^n} \rho_1(\Yb_n) \left( \R \backslash \left\{ \int_0^1 f_0(x)\dd x \right\} \right) \Pb_{1,n}^{f_0}(\dd \Yb_n)\\
&\geq \inf_{\rho_1} \int_{A_0} \rho_1(\Yb_n) \left( \R \backslash \left\{ \int_0^1 f_0(x)\dd x \right\} \right) \Pb_{1,n}^{f_0}(\dd \Yb_n)\\
&\geq \frac{1}{2} \Pb_{1,n}^{f_0}(A_0)\\
&\geq \frac{1}{4},
\end{align*}
which proves the assertion.   
\printbibliography

\end{document}